\documentstyle[12pt]{article}
\input amssym.def
\begin{document}
\def \Z{\Bbb Z}
\def \C{\Bbb C}
\def \R{\Bbb R}
\def \Q{\Bbb Q}
\def \N{\Bbb N}
\def \P{\Bbb P}
\def \bZ{{\bf Z}}
\def \TZ{{1\over T}{\bf Z}}
\def \wt{{\rm wt}\;}
\def \gr{{\rm gr}}
\def \g{\frak{g}}
\def \gl{\frak{gl}}
\def \h{\frak h}
\def \mod{{\rm mod}\;}
\def \pf{{\bf Proof}.}

\def \Der{{\rm Der}}
\def \Sing{{\rm Sing}}
\def \span{{\rm span}}
\def \Res{{\rm Res}}
\def \End{{\rm End}}
\def \Hom{{\rm Hom}}
\def \<{\langle} 
\def \>{\rangle}
\def \be{\begin{equation}\label}
\def \ee{\end{equation}}
\def \bex{\begin{exa}\label}
\def \eex{\end{exa}}
\def \bl{\begin{lem}\label}
\def \el{\end{lem}}
\def \bt{\begin{thm}\label}
\def \et{\end{thm}}
\def \bp{\begin{prop}\label}
\def \ep{\end{prop}}
\def \br{\begin{rem}\label}
\def \er{\end{rem}}
\def \bc{\begin{coro}\label}
\def \ec{\end{coro}}
\def \bd{\begin{de}\label}
\def \ed{\end{de}}

\newtheorem{thm}{Theorem}[section]
\newtheorem{prop}[thm]{Proposition}
\newtheorem{coro}[thm]{Corollary}
\newtheorem{conj}[thm]{Conjecture}
\newtheorem{exa}[thm]{Example}
\newtheorem{lem}[thm]{Lemma}
\newtheorem{rem}[thm]{Remark}
\newtheorem{de}[thm]{Definition}
\newtheorem{hy}[thm]{Hypothesis}
\makeatletter
\@addtoreset{equation}{section}
\def\theequation{\thesection.\arabic{equation}}
\makeatother
\makeatletter

\begin{center}
{\Large \bf Vertex algebras and vertex Poisson algebras}
\end{center}
\begin{center}{Haisheng Li\footnote{Partially supported by NSF grant
DMS-9970496 and a grant from Rutgers Research Council}\\
Department of Mathematical Sciences, Rutgers University, Camden, NJ 08102\\
and\\
Department of Mathematics, Harbin Normal University, Harbin, China}
\end{center}

\begin{abstract}
This paper studies certain relations among
vertex algebras, vertex Lie algebras
and vertex Poisson algebras.
In this paper, the notions of vertex Lie algebra (conformal algebra)
and vertex Poisson algebra are revisited and
certain general construction theorems
of vertex Poisson algebras are given.
A notion of filtered vertex algebra is formulated 
in terms of a notion of good filtration and 
it is proved that the associated graded vector space
of a filtered vertex algebra 
is naturally a vertex Poisson algebra.
For any vertex algebra $V$, a general construction and 
a classification of good filtrations are given.
To each $\N$-graded vertex algebra $V=\coprod_{n\in \N}V_{(n)}$ 
with $V_{(0)}=\C {\bf 1}$, a canonical
(good) filtration is associated and certain results 
about generating subspaces of certain types of $V$ 
are also obtained.
Furthermore, a notion of formal deformation of a 
vertex (Poisson) algebra is formulated and 
a formal deformation of vertex Poisson algebras
associated with vertex Lie algebras is constructed.
\end{abstract}

\section{Introduction}
In quantum physics, a very important role is played by
Poisson algebras. A Poisson algebra structure 
on a vector space $A$ is a combination of 
a commutative associative algebra structure
and a Lie algebra structure with a certain 
compatibility condition (the Leibniz rule).
Important examples from quantum mechanics are those with $A=C^{\infty}(M)$, 
where $M$ are symplectic manifolds (the phase spaces of dynamical systems).
Another important family of Poisson algebras 
are the Poisson algebras $S(\g)$ $(=\C[\g^{*}])$ associated with 
(finite-dimensional) Lie algebras $\g$.
In the Poisson world, among the most important issues is 
the algebraic deformation quantization of Poisson algebras
(see [BFFLS], [Ko]).

Vertex (operator) algebras are known as the mathematical
counterparts of chiral algebras in two-dimensional quantum
conformal field theory and are analogous to associative algebras
in certain aspects. (At the same time, in many
aspects vertex (operator) algebras are also
analogous to both Lie algebras and commutative associative algebras.) 
One example of such analogy 
is that Borcherds' commutator formula, 
which is just {\em part} of the vertex algebra structure, 
gives rise to a certain Lie algebra (cf. [B1], [FFR], [Li1], [MP]),
just as the commutator of an associative algebra 
gives rise to a canonical Lie algebra.
This essentially motivated the introduction of
Lie algebra analogues of vertex algebras.
In [K], among other things Kac introduced a notion of conformal algebra
and independently Primc introduced and studied
a notion of vertex Lie algebra in [P].
(A different notion of vertex Lie algebra was 
introduced in [DLM2].) As it was explained in 
Remark \ref{requivalence-vla-conformala},
the notion of conformal algebra and the notion of vertex Lie algebra
are equivalent. 
Following [FB] we use conformal algebra and vertex Lie algebra
synonymously in this paper. 
The main axiom defining the notion of vertex Lie
algebra is what Primc called ``half commutator formula.''
Associated to each vertex Lie algebra $R$,
there is an honest Lie algebra ${\cal{L}}(R)$, the underlying vector
space of which is a certain quotient space of the loop space $L(R)$ 
$(=R\otimes \C[t,t^{-1}])$ of $R$.
Just as one gets the universal enveloping
algebra $U(\g)$ (an associative algebra) from any Lie algebra $\g$, 
one gets a canonical vertex algebra
${\cal{V}}(R)$ (see [DLM2], [FB], [K], [P]) from 
any vertex Lie algebra $R$, where the vertex algebra ${\cal{V}}(R)$
is a suitably defined Verma ${\cal{L}}(R)$-module.

With the notion of vertex Lie algebra, 
one naturally arrives at the notion of vertex Poisson algebra.
A vertex Poisson algebra structure (see [FB]; cf. [BD], [EF], [DLM2])
is a combination of a commutative vertex algebra structure (or
equivalently a differential algebra structure)
and a vertex Lie algebra structure with a natural compatibility condition.
Vertex Poisson algebras were studied in [FB] (Chapter 15)
and many important results were obtained.
Among those results, 
it was noticed that the symmetric algebra of
any vertex Lie algebra is naturally a vertex Poisson algebra and 
certain vertex Poisson algebras associated with
affine and the Virasoro Lie algebras were 
realized as classical limits of vertex algebras.
Certain connections between the classical and 
the quantum Drinfeld-Sokolov reductions were also exhibited.

In this paper, we shall study the connection of vertex algebras,
vertex Lie algebras and vertex Poisson algebras 
in a more systematic way and the main purpose of this paper
is to lay the foundation for future studies 
on vertex algebras and vertex Poisson algebras.
As one of our main results, we introduce and study a notion of 
filtered vertex algebra and we show that
the graded vector space of a filtered vertex algebra 
is naturally a vertex Poisson algebra. For any vertex algebra $V$,
we give a general construction and a classification of filtered vertex
algebras $(V,E)$ and we associate a canonical filtered vertex algebra
to every $\N$-graded vertex algebra $V$ with $V_{(0)}=\C {\bf 1}$.
The introduction of the notion of filtered vertex algebra was motivated 
by certain classical results and certain results in [KL]
on generating subspaces of a certain type of vertex operator algebras.
In the classical case, the standard filtration on the universal
enveloping algebra of a Lie algebra and 
the fundamental Poincare-Birkhoff-Witt theorem
are closely related. Similarly, the canonical filtration we constructed
is closely related to what we call generating subspaces with PBW
spanning property and certain basic results are obtained.

For determining vertex Poisson algebra structures on a given differential
algebra we give two general construction theorems of vertex Poisson algebras.
Such construction theorems are very effective to determine
whether a certain vertex Poisson algebra structure exists on a 
differential algebra. In particular, 
our construction theorems can be applied to show 
that there exists a natural vertex Poisson
algebra structure on the symmetric algebra $S(R)$ of
a vertex Lie algebra $R$ (Proposition \ref{panother}; 
see also [FB], Example 15.2.2).
As it is needed in the proof of the first construction theorem,
we prove that in a suitable setting, the half commutator formula
has $S_{3}$-symmetry. 
This makes vertex Lie algebras virtually look like Lie algebras.
The proof of this particular result
is based on certain results of [P] and [FHL].

For any vertex Lie algebra $R$, we construct a canonical
map from the symmetric
algebra $S(R)$ onto ${\cal{V}}(R)$, which is analogous to
the canonical symmetrization map from $S(\g)$ to $U(\g)$ 
for a Lie algebra $\g$, and we show that the vertex Poisson algebra
$S(R)$ is isomorphic to $\gr_{E}{\cal{V}}(R)$, 
where $E$ is the associated canonical filtration.
Motivated by a result of Frenkel and Ben-Zvi [FB], 
we formulate a notion of $*$-deformation of a vertex Poisson
algebra and we construct a $*$-deformation of
the vertex Poisson algebras associated with filtered vertex algebras
and with vertex Lie algebras.

As with Poisson algebras, a fundamental problem is
about the existence, uniqueness and construction of
each vertex Poisson algebra.
In a sequel we shall study the existence
of $*$-deformation of a vertex Poisson structure 
on a free differential algebra
(a polynomial algebra) and we shall 
use the vertex Poisson algebra $\gr_{E}\; V$ 
to study the algebraic structures of a vertex algebra $V$. 

In this work we use many of Primc's results
and ideas from [P] and we certainly use many ideas from 
[DLM2]. We benefit very much from 
reading Chapter 15 of [FB] and from talking with Ping Xu.

This paper is organized as follows: In Section 2,
we review and discuss the notion of vertex Lie algebra 
(conformal algebra) and we present some new results.
In Section 3, we recall the notion of vertex Poisson algebra 
and give two general construction theorems of vertex Poisson algebras.
In Section 4, we study filtered vertex algebras.
In Section 5, we study $*$-deformations of
vertex (Poisson) algebras.

\section{Vertex algebras and vertex Lie algebras} 
In this section we review and discuss the notions of vertex
algebra and vertex Lie algebra (conformal algebra), and their relations.
We recall certain results of Primc about the half commutator
formula and the half Jacobi identity and we prove that
the half Jacobi identity and the half commutator
formula have $S_{3}$-symmetry.
We also prove a result analogous to the classical result about 
the canonical (symmetrization) map from $S(\g)$ to $U(\g)$
for a Lie algebra $\g$.

In this paper we use the standard formal variable notations 
and conventions as defined in [FLM] and [FHL] (cf. [LL]).
In addition to the standard notations $\C$ and $\Z$,
 we use notations $\N$ 
for the set of nonnegative integers
and $\Z_{+}$ for the set of positive integers.

Let $V$ be any vector space. Following [P], 
for a formal series
$$f(x_{1},\dots,x_{n})=\sum_{m_{1},\dots,m_{n}\in \Z}u(m_{1},\dots,m_{n})
x_{1}^{-m_{1}-1}\cdots x_{n}^{-m_{n}-1}
\in V[[x_{1}^{\pm 1},\dots,x_{n}^{\pm 1}]],$$
we set
\begin{eqnarray}
\Sing f(x_{1},\dots,x_{n})
=\sum_{m_{1},\dots,m_{n}\ge 0}u(m_{1},\dots,m_{n})
x_{1}^{-m_{1}-1}\cdots x_{n}^{-m_{n}-1}.
\end{eqnarray}
Clearly, for $1\le i\le n$, 
\begin{eqnarray}
{\partial\over\partial x_{i}}\Sing f(x_{1},\dots,x_{n})
=\Sing {\partial\over\partial x_{i}}f(x_{1},\dots,x_{n}).
\end{eqnarray}
For any nonempty subset $S=\{ i_{1},\dots, i_{k}\}$ of $\{ 1,\dots, n\}$, 
a formal series $f(x_{1},\dots,x_{n})$ can be naturally viewed
as a formal series $\tilde{f}(x_{i_{1}},\dots,x_{i_{k}})$ 
in variables $x_{i_{1}},\dots, x_{i_{k}}$ with coefficients in
the vector space $V[[x_{j}^{\pm 1}\;|\; j\notin S]]$. Then we define
\begin{eqnarray}
\Sing_{x_{i_{1}},\dots,x_{i_{k}}} f(x_{1},\dots,x_{n})
=\Sing \tilde{f}(x_{i_{1}},\dots,x_{i_{k}}).
\end{eqnarray}
In particular, we have
\begin{eqnarray}
\Sing_{x_{1}} f(x_{1},\dots,x_{n})
=\sum_{m\in \N,\; m_{2},\dots,m_{n}\in \Z}u(m,m_{2},\dots,m_{n})
x_{1}^{-m-1}x_{2}^{-m_{2}-1}\cdots x_{n}^{-m_{n}-1}.
\end{eqnarray}
Then
\begin{eqnarray}
\Sing f(x_{1},\dots,x_{n})
=\Sing_{x_{1}}\cdots \Sing_{x_{n}}f(x_{1},\dots,x_{n}).
\end{eqnarray}
The following result is due to Primc in [P]:

\bl{lsingularformula}
Let $V$ be a vector space and let 
$$B\in V((x_{1},\dots,x_{n})),\;
P\in (\End \; V)[[x_{1},\dots,x_{n}]].$$ 
Then
\begin{eqnarray}
\Sing \left(P\cdot \Sing (B)\right)=\Sing(PB).
\end{eqnarray}
\el

In this paper we use the following definition of the notion of 
vertex algebra as in [LL]:

\bd{dvadefinition}
{\em A {\em vertex algebra} is a vector space $V$ equipped with 
a linear map 
\begin{eqnarray}
Y: V&\rightarrow &\Hom (V,V((x)))\nonumber\\
v&\mapsto &Y(v,x)=\sum_{n\in \Z}v_{n}x^{-n-1}\;\;(\mbox{\rm where}\;
v_{n}\in \End V)
\end{eqnarray}
and equipped with a distinguished vector ${\bf 1}\in V$,
called the {\em vacuum vector},
such that the following axioms hold: 
\begin{eqnarray}\label{evacuumva}
Y({\bf 1},x)=1;
\end{eqnarray}
for $v\in V$, 
\begin{eqnarray}
Y(v,x){\bf 1}\in V[[x]]\;\;\;\mbox{ and }\;\;\; 
Y(v,x){\bf 1}|_{x=0}\left(=v_{-1}{\bf 1}\right)=v;
\end{eqnarray}
and for $u,v\in V$,
\begin{eqnarray}\label{ejacobiva}
& &x_{0}^{-1}\delta\left(\frac{x_{1}-x_{2}}{x_{0}}\right)
Y(u,x_{1})Y(v,x_{2})-x_{0}^{-1}\delta\left(\frac{x_{2}-x_{1}}{-x_{0}}\right)
Y(v,x_{2})Y(u,x_{1})\nonumber\\
&=&x_{2}^{-1}\delta\left(\frac{x_{1}-x_{0}}{x_{2}}\right)
Y(Y(u,x_{0})v,x_{2})
\end{eqnarray}
(the {\em Jacobi identity}).}
\ed

By taking $\Res_{x_{0}}$ and $\Res_{x_{1}}$ respectively from 
the Jacobi identity 
we get Borcherds' commutator formula and iterate formula
\begin{eqnarray}
& &[Y(u,x_{1}),Y(v,x_{2})]
=\Res_{x_{0}}x_{2}^{-1}\delta\left(\frac{x_{1}-x_{0}}{x_{2}}\right)
Y(Y(u,x_{0})v,x_{2}),\\
& &Y(Y(u,x_{0})v,x_{2})\nonumber\\
&=&\Res_{x_{1}}\left(x_{0}^{-1}\delta\left(\frac{x_{1}-x_{2}}{x_{0}}\right)
Y(u,x_{1})Y(v,x_{2})-x_{0}^{-1}\delta\left(\frac{x_{2}-x_{1}}{-x_{0}}\right)
Y(v,x_{2})Y(u,x_{1})\right).\;\;\;\;\;
\end{eqnarray}
In terms of components, we have 
\begin{eqnarray}
& &[u_{m},v_{n}]=\sum_{i\ge 0}{m\choose i}(u_{i}v)_{m+n-i},
\label{ecommutatorincomponents}\\
& &(u_{m}v)_{n}=\sum_{i\ge 0}{m\choose i}(-1)^{i}
\left( u_{m-i}v_{n+i}-(-1)^{m}v_{m+n-i}u_{i}\right)
\label{eiterateincomponents}
\end{eqnarray}
for $m,n\in \Z$. 
Let ${\cal{D}}$ be the linear operator on $V$ defined by
\begin{eqnarray}\label{edoperator-definition}
{\cal{D}}(v)=v_{-2}{\bf 1}
\left(=\left({d\over dx}Y(v,x){\bf 1}\right)|_{x=0}\right).
\end{eqnarray}
Then (cf. [LL])
\begin{eqnarray}
& &[{\cal{D}},Y(u,x)]=Y({\cal{D}}u,x)={d\over dx}Y(u,x),\\
& &Y(u,x)v=e^{x{\cal{D}}}Y(v,-x)u\;\;\;\mbox{ for }u,v\in V.
\end{eqnarray}
In terms of components, we have
\begin{eqnarray}
& &[{\cal{D}},u_{m}]=({\cal{D}}u)_{m}=-m u_{m-1}
\label{edbracketformulaincomponenets}\\
& & u_{m}v=\sum_{i\ge 0}(-1)^{m+i-1}{1\over i!}{\cal{D}}^{i} v_{m+i}u
\label{eskew-symmetry-components}
\end{eqnarray}
for $u,v\in V,\; m\in \Z$.

It follows from Borcherds' commutator formula 
(\ref{ecommutatorincomponents}) that all the operators $u_{m}$ on $V$,
 for $u\in V,\; m\in \Z$, linearly span a Lie subalgebra $\g_{V}$ of
$\gl(V)$ (the general linear Lie algebra).
On the other hand, form the vector space
\begin{eqnarray}
L(V)=V\otimes \C[t,t^{-1}]
\end{eqnarray}
and define a bilinear multiplication on $L(V)$ by
\begin{eqnarray}
[u\otimes t^{m},v\otimes t^{n}]
=\sum_{i\ge 0} {m\choose i}(u_{i}v\otimes t^{m+n-i})
\end{eqnarray}
for $u,v\in V,\; m,n\in \Z$.
It was proved (see [FFR], [Li1], [MP]) that
the multiplication on $L(V)$ 
reduces to a Lie algebra structure on the quotient space
\begin{eqnarray}
{\cal{L}}(V)=L(V)/DL(V),
\end{eqnarray}
where $D={\cal{D}}\otimes 1+1\otimes {d\over dt}$.
The Lie algebra $\g_{V}$ is naturally a homomorphic image of ${\cal{L}}(V)$.
Notice that the Lie algebra structures on $\g_{V}$ and ${\cal{L}}(V)$ 
only use the singular part of $Y(u,x)v$ for $u,v\in V$. 
This essentially motivated the introduction of
the notions of conformal algebra and vertex Lie algebra.
Following [FB] {\em we use conformal algebra and vertex Lie algebra
synonymously in this paper}.

The following notion of vertex Lie algebra 
is due to [K] under the name of
conformal algebra and independently due to [P] 
(see Remark \ref{requivalence-vla-conformala}):

\bd{dkacconformalalgebra}
{\em A {\em vertex Lie algebra} is a vector space $A$ equipped with
a linear operator $\partial$ and equipped with $\C$-bilinear products
$(a,b)\mapsto a_{n}b$ for $a,b\in A,\; n\in \N$, such that 
the following axioms hold for $a,b,c\in A,\; m,n\in \N$:

(C0) $a_{n}b=0$ for $n$ sufficiently large;

(C1) $(\partial a)_{n}b=-na_{n-1}b$;

(C2) 
\begin{eqnarray}
a_{n}b=\sum_{i\ge 0}(-1)^{n+i+1}{1\over i!}\partial^{i}b_{n+i}a;
\end{eqnarray}

(C3)
\begin{eqnarray}
a_{m}b_{n}c-b_{n}a_{m}c=\sum_{i=0}^{m}{m\choose i}(a_{i}b)_{m+n-i}c.
\end{eqnarray}}
\ed

\br{rva-ca}
{\em In view of (\ref{edbracketformulaincomponenets}), 
(\ref{eskew-symmetry-components}) and (\ref{ecommutatorincomponents}),
any vertex algebra $V$ is naturally a 
vertex Lie algebra with $\partial ={\cal{D}}$, the operator defined in
(\ref{edoperator-definition}).}
\er

In the notion of vertex Lie algebra, the $\C$-bilinear products 
(in terms of generating functions) and the axiom (C0) 
amount to a linear map 
\begin{eqnarray*}
Y_{-}: A&\rightarrow& \Hom (A,x^{-1}A[x^{-1}])\\
a&\mapsto &Y_{-}(a,x)=\sum_{n\ge 0}a_{n}x^{-n-1}\;\;\;\mbox{(where
$a_{n}\in \End\; A$)}.
\end{eqnarray*}
With the axiom (C0) being included in
the map $Y_{-}$, the axioms (C1)-(C3) amount to
\begin{eqnarray}
& &Y_{-}(\partial a,x)={d\over dx}Y_{-}(a,x)\label{ederivativeproperty}\\
& &Y_{-}(a,x)b=\Sing\left(e^{x\partial}Y_{-}(b,-x)a\right)
\label{eskewsymmetry}\\
& &[Y_{-}(a,x_{1}),Y_{-}(b,x_{2})]
=\Sing\left(\sum_{i\ge 0}(x_{1}-x_{2})^{-i-1}Y_{-}(a_{i}b,x_{2})\right)\\
& &\hspace{3.75cm}=\Sing\left(Y_{-}(Y_{-}(a,x_{1}-x_{2})b,x_{2})\right).
\label{ehalfcommutatorformula1}
\end{eqnarray}

Differentiating both sides of (\ref{eskewsymmetry}) and then using 
the properties (\ref{eskewsymmetry}) and (\ref{ederivativeproperty}) we get
\begin{eqnarray*}
{d\over dx}Y_{-}(a,x)b
&=&\Sing \left(\partial e^{x\partial}Y_{-}(b,-x)a-e^{x\partial}
Y_{-}(\partial b,-x)a\right)\nonumber\\
&=&\partial Y_{-}(a,x)b-Y_{-}(a,x)\partial b.
\end{eqnarray*}
Thus (see [K])
\begin{eqnarray}\label{ethreetogether}
[\partial,Y_{-}(a,x)]={d\over dx}Y_{-}(a,x)=Y_{-}(\partial a,x)
\;\;\;\mbox{ for }a\in A.
\end{eqnarray}
Equivalently we have 
\begin{eqnarray}\label{econjugation-formula}
e^{x_{1}\partial}Y_{-}(a,x)e^{-x_{1}\partial}
=e^{x_{1}{d\over dx}}Y_{-}(a,x)=Y_{-}(e^{x_{1}\partial} a,x)=Y_{-}(a,x+x_{1})
\end{eqnarray}
for $a\in A$. Following [P], we call (\ref{ehalfcommutatorformula1}) the 
{\em half commutator formula} and call the following identity
the {\em half Jacobi identity}:
\begin{eqnarray}\label{ehalfjacobiidentity}
& &\Sing\left(x_{0}^{-1}\delta\left(\frac{x_{1}-x_{2}}{x_{0}}\right)
Y_{-}(a,x_{1})Y_{-}(b,x_{2})
-x_{0}^{-1}\delta\left(\frac{x_{2}-x_{1}}{-x_{0}}\right)
Y_{-}(b,x_{2})Y_{-}(a,x_{1})\right)\nonumber\\
&=& \Sing\left( x_{2}^{-1}\delta\left(\frac{x_{1}-x_{0}}{x_{2}}\right)
Y_{-}(Y_{-}(a,x_{0})b,x_{2})\right).
\end{eqnarray}
We refer to (\ref{ehalfjacobiidentity}) applied to a vector $c\in A$
as the {\em half Jacobi identity for the triple $(a,b,c)$}
and similarly we refer to (\ref{ehalfcommutatorformula1})
applied to a vector $c\in A$
as the {\em half commutator formula for the triple $(a,b,c)$}.

The following result is due to Primc ([P], Lemma 6.1;
notice that the $D$-operator in the setting of Lemma 6.1 was not used
in the proof):

\bp{pprimcequivalence}
Let $A$ be a vector space equipped with a linear map $Y_{-}$ from $A$ to 
$\Hom (A,x^{-1}A[x^{-1}])$. Then the half commutator relation 
(\ref{ehalfcommutatorformula1}) is equivalent
to the half Jacobi identity relation (\ref{ehalfjacobiidentity}).
\ep

\br{requivalence-vla-conformala}
{\em In [P], a vertex Lie algebra was defined 
to be a vector space $A$ equipped with a linear
operator $D$ and equipped with a linear map $Y_{-}$ from $A$ to $\Hom
(A,x^{-1}A[x^{-1}])$ such that (\ref{eskewsymmetry}), (\ref{ethreetogether})
with $D$ in place of $\partial$ and (\ref{ehalfjacobiidentity}) hold.
In view of Proposition \ref{pprimcequivalence}, this definition 
is equivalent to Definition \ref{dkacconformalalgebra}.
We often denote a vertex Lie algebra
by $(A, Y_{-},\partial)$ and refer to $(Y_{-},\partial)$ 
as the vertex Lie algebra structure.}
\er

Now we discuss a certain symmetry
of the half skew symmetry and the half Jacobi identity,
which will be useful in constructing vertex Lie (Poisson) algebras.

\bl{lskewsymmetry1}
Let $V$ be a vector space equipped with a linear operator $\partial$
and let $A(x), B(x)\in x^{-1}V[x^{-1}]$.
Then
\begin{eqnarray}\label{eskewsymmetryab}
A(x)=\Sing\left(e^{x\partial}B(-x)\right)
\end{eqnarray}
if and only if
\begin{eqnarray}\label{eskewsymmetryba}
B(x)=\Sing\left(e^{x\partial}A(-x)\right).
\end{eqnarray}
\el

{\bf Proof.} From the obvious symmetry we only need to prove that
(\ref{eskewsymmetryab}) implies (\ref{eskewsymmetryba}). 
Using (\ref{eskewsymmetryab}) and Lemma \ref{lsingularformula} we obtain
(\ref{eskewsymmetryba}) as
\begin{eqnarray}
\Sing\left(e^{x\partial}A(-x)\right)
=\Sing\left(e^{x\partial}\Sing\left(e^{-x\partial}B(x)\right)\right)
=\Sing B(x)=B(x).\;\;\;\;\Box
\end{eqnarray}

The following is a simple observation from the $S_{3}$-symmetry
of the Jacobi identity for vertex algebras (see [FHL]): 

\bp{ps3symmetry}
Let $A$ be a vector space equipped with a linear operator $\partial$ 
and equipped with a linear map 
\begin{eqnarray}
Y_{-}: A&\rightarrow& \Hom (A,x^{-1}A[x^{-1}])\nonumber\\
a&\mapsto& Y_{-}(a,x)
\end{eqnarray}
such that the following conditions hold for all $a,b\in A$:
\begin{eqnarray}
& &[\partial,Y_{-}(a,x)]={d\over dx}Y_{-}(a,x),\label{ebracketderivative}\\
& &Y_{-}(a,x)b=\Sing\left(e^{x\partial}Y_{-}(b,-x)a\right).
\end{eqnarray}
Then the half Jacobi identity for an ordered triple $(a,b,c)$
implies the half Jacobi identity for any permutation of $(a,b,c)$.
\ep

{\bf Proof.} Notice that from (\ref{ebracketderivative}) we 
have the following conjugation formula
\begin{eqnarray}\label{econjugation}
e^{x_{1}\partial}Y_{-}(a,x)e^{-x_{1}\partial}
=Y_{-}(a, x+x_{1})\;\;\;\mbox{ for }a\in A.
\end{eqnarray}
Then using the same argument of [FHL] (and using Lemma \ref{lsingularformula},
the half skew symmetry and the conjugation formula) we have
the $S_{3}$-symmetry of the half Jacobi identity. $\;\;\;\;\Box$

Combining Propositions \ref{pprimcequivalence} and \ref{ps3symmetry}
we immediately have:

\bp{ps3symmetry-commutator}
Let $A$, $\partial$ and $Y_{-}$ be as in Proposition \ref{ps3symmetry}
and let $a,b,c\in A$.
Then the half commutator formula for $(a,b,c)$ implies 
the half commutator formula for any permutation of $(a,b,c)$.$\;\;\;\;\Box$
\ep

The following result, which is due to Primc [P],
relates a vertex Lie algebra canonically with
an honest Lie algebra:

\bp{pliealgebra}
Let $A$ be a vector space equipped with a linear operator $\partial$ and 
equipped with a linear map $Y_{-}$ from $A$ to $\Hom (A,x^{-1}A[x^{-1}])$,
such that
\begin{eqnarray}
[\partial,Y_{-}(a,x)]=Y_{-}(\partial a,x)={d\over dx}Y_{-}(a,x)
\;\;\;\mbox{ for }a\in A.
\end{eqnarray}
Set
\begin{eqnarray}
L(A)=A\otimes \C[t,t^{-1}]
\end{eqnarray}
and set
\begin{eqnarray}
\hat{\partial}=\partial\otimes 1+1\otimes {d\over dt}\in \End\; L(A).
\end{eqnarray}
Define
\begin{eqnarray}
[a\otimes t^{m},b\otimes t^{n}]
=\sum_{i\ge 0}{m\choose i}a_{i}b\otimes t^{m+n-i}
\end{eqnarray}
for $a,b\in A,\; m,n\in \N$, where $Y_{-}(a,x)b=\sum_{n\ge 0}a_{n}b x^{-n-1}$.
Then $A$ is a vertex Lie algebra 
if and only if $L(A)/\hat{\partial}L(A)$ is a Lie algebra.
\ep

Let $R$ be a vertex Lie algebra.
Set 
\begin{eqnarray}
{\cal{L}}(R)=L(R)/\hat{\partial}L(R),
\end{eqnarray}
the Lie algebra associated to $R$ through Proposition \ref{pliealgebra}.
Denote by $\rho$ the quotient map from $L(R)$ to 
${\cal{L}}(R)$. Then $\ker \rho =\hat{\partial}L(R)$.

\br{rdlm-definition}
{\em A different notion of vertex Lie algebra was formulated in [DLM2]
with the Lie algebra structure and the surjective map
$\rho$ as the main components.}
\er

\br{rcategory-Liealg}
{\em It is straightforward to see that a vertex Lie algebra 
homomorphism from $R_{1}$ to $R_{2}$ naturally
gives rise to a Lie algebra homomorphism from ${\cal{L}}(R_{1})$ 
to ${\cal{L}}(R_{2})$. Then we have a functor from the category of
vertex Lie algebras to the category of Lie algebras.}
\er

Let $R$ be a vertex Lie algebra. 
Then we have a canonical polar decomposition of the associated Lie
algebra (see [DLM2], [P]):
\begin{eqnarray}
{\cal{L}}(R)={\cal{L}}(R)_{+}\oplus {\cal{L}}(R)_{-},
\end{eqnarray}
where 
\begin{eqnarray}
{\cal{L}}(R)_{-}=\rho\left(R\otimes t^{-1}\C[t^{-1}]\right)
\;\;\;\mbox{ and }\;\;\; 
{\cal{L}}(R)_{+}=\rho\left(R\otimes \C[t]\right)
\end{eqnarray}
are subalgebras. For $a\in R,\; n\in \Z$, set
\begin{eqnarray}
a(n)=\rho(a\otimes t^{n})\in {\cal{L}}(R).
\end{eqnarray}

\br{rhalfmodulestructure}
{\em For any ${\cal{L}}(R)$-module $W$, we also use
$a(n)$ for the corresponding operator on $W$.
It is clear that $R$ is a natural ${\cal{L}}(R)_{+}$-module
with $a(n)$ acting as $a_{n}$ for $a\in R,\; n\ge 0$.}
\er

\br{rborcherds}
{\em Let $R$ be a vertex Lie algebra.
Then $R/\partial R$ is a Lie algebra
with 
\begin{eqnarray}
[a+\partial R,b+\partial R]=a_{0}b+\partial R
\end{eqnarray}
for $a,b\in R$ (see [B1]). On the other hand, it is easy to see that
$\rho (R)$ is a Lie subalgebra of ${\cal{L}}(R)_{+}$ and $\ker \rho \cap
R=\partial R$. Then $R/\partial R$ is isomorphic to
the Lie subalgebra $\rho (R)$ of ${\cal{L}}(R)_{+}$.
We denote the Lie algebra $R/\partial R=\rho (R)$ by ${\cal{L}}(R)_{0}$.}
\er

An {\em $R$-module} (see [K]) is an ${\cal{L}}(R)_{+}$-module $W$ 
with an action of $\partial$ such that
\begin{eqnarray}
[\partial, Y_{-}(u,x)]={d\over dx}Y_{-}(\partial u,x)
\;\;\;\mbox{ for }u\in R.
\end{eqnarray}

View $\C$ as a trivial ${\cal{L}}(R)_{+}$-module and form the induced module
\begin{eqnarray}
{\cal{V}}(R)=U({\cal{L}}(R))\otimes _{U({\cal{L}}(R)_{+}}\C.
\end{eqnarray}
In view of the Poincar\'e-Birkhoff-Witt theorem, we have
\begin{eqnarray}
{\cal{V}}(R)=U({\cal{L}}(R)_{-}),
\end{eqnarray}
as a vector space, so that
we may and do consider $\C$ as a subspace of ${\cal{V}}(R)$. Set
\begin{eqnarray}
{\bf 1}=1\otimes 1\in {\cal{V}}(R).
\end{eqnarray}
For the same reason, we consider $R$ as a subspace 
of ${\cal{V}}(R)$ through the map $a\mapsto a(-1){\bf 1}$.

The following result was proved in [DLM2] and [P] (see also [FB] and
[K]):

\bt{tdlmkp}
There exists a unique vertex algebra structure
on the ${\cal{L}}(R)$-module ${\cal{V}}(R)$ such that
${\bf 1}$ is the vacuum vector and such that
\begin{eqnarray}
Y(a,x)=\sum_{n\in \Z}a(n)x^{-n-1}\;\;\;\mbox{ for }a\in R.
\end{eqnarray}
The subspace $R$ generates ${\cal{L}}(R)$ as a vertex algebra.
Furthermore, any restricted ${\cal{L}}(R)$-module $W$,
in the sense that for any $w\in W$ and $a\in R$, $a(n)w=0$ 
for $n$ sufficiently large, is a natural ${\cal{V}}(R)$-module.
\et

The following result is due to Primc in [P]:

\bt{tprimc}
Let $R$ be a vertex Lie algebra and let ${\cal{V}}(R)$ be
the vertex algebra associated with $R$. Then the identification map
of $R$ as a subspace of ${\cal{V}}(R)$ is a vertex Lie algebra
homomorphism. Furthermore,
for any vertex algebra $V$ and for any vertex Lie algebra 
homomorphism $f$ from $R$ to $V$ viewed as a vertex Lie algebra,
there exists a unique vertex algebra homomorphism $\bar{f}$ from
${\cal{V}}(R)$ to $V$, extending $f$.
\et

The following result is due to Primc [P] (it was also proved in [DLM2] for
$R=V$ a vertex algebra):

\bp{pprimc}
Let $R$ be a vertex Lie algebra. Then the linear map
\begin{eqnarray}
R&\rightarrow& {\cal{L}}(R)_{-}\nonumber\\
a&\mapsto& a(-1)
\end{eqnarray}
is a bijection. 
\ep

Since $R$ is an ${\cal{L}}(R)_{+}$-module, the symmetric algebra $S(R)$ 
is naturally an ${\cal{L}}(R)_{+}$-module. 
On the other hand, ${\cal{V}}(R)$ is also
an ${\cal{L}}(R)_{+}$-module. In particular, both $S(R)$ 
and ${\cal{V}}(R)$ are natural ${\cal{L}}(R)_{0}$-modules.
The following is an analogue of
a classical result in Lie theory (cf. [Di]):

\bp{psymmetrizationmap}
Let $R$ be a vertex Lie algebra.
Then the linear map $\omega$ from $S(R)$ to ${\cal{V}}(R)$, defined by
\begin{eqnarray}
\omega (u^{(1)}\cdots u^{(r)})=
{1\over r!}\sum_{\sigma\in S_{r}}u^{(\sigma(1))}(-1)\cdots
u^{(\sigma(r))}(-1){\bf 1}
\end{eqnarray}
for $r\ge 0,\;u^{(i)}\in R$, is an ${\cal{L}}(R)_{0}$-module
isomorphism,
where $S_{r}$ is the symmetric group on $\{1,\dots,r\}$.
\ep

{\bf Proof.} In view of Proposition \ref{pprimc}
and the Poincar\'e-Birkhoff-Witt theorem, we have
a linear isomorphism 
\begin{eqnarray}
f: U({\cal{L}}(R)_{-})&\rightarrow& {\cal{V}}(R)\nonumber\\
   u^{(1)}(-1)\cdots u^{(r)}(-1)&\mapsto& u^{(1)}(-1)\cdots u^{(r)}(-1){\bf 1}
\end{eqnarray}
for $r\ge 0$, $u^{(i)}\in R$. View $R$ as a Lie algebra
with the transported structure from ${\cal{L}}(R)_{-}$.
{}From [Di], the canonical (symmetrization) map from $S(R)$ to $U(R)$
is a linear isomorphism. Consequently, 
the linear map $\omega$ is a linear isomorphism.

Furthermore, for $u,v\in R,\; w\in {\cal{V}}(R)$, from Borcherds' commutator
formula we have
\begin{eqnarray}\label{evmu-1w}
v(0)u(-1)w=u(-1)v(0)w+(v_{0}u)(-1)w.
\end{eqnarray}
Since $v(0){\bf 1}=0$, it follows that $\omega$ is 
an ${\cal{L}}(R)_{0}$-module homomorphism.
Therefore, $\omega$ is an ${\cal{L}}(R)_{0}$-module isomorphism.
$\;\;\;\;\Box$

For a module $W$ for a Lie algebra $\g$, we set
\begin{eqnarray}
W^{\g}=\{w\in W\;|\; \g w=0\}.
\end{eqnarray}

As an immediate consequence of Proposition \ref{psymmetrizationmap}
we have:

\bc{cwalgebra1}
Let $\g$ be a Lie subalgebra of ${\cal{L}}(R)_{0}$. 
Then the linear map $\omega$ defined in Proposition \ref{psymmetrizationmap}
gives rise to a linear isomorphism from
$S(R)^{\g}$ onto ${\cal{V}}(R)^{\g}$.$\;\;\;\;\Box$
\ec

\br{rwalgebra2}
{\em A {\em derivation} (see [Lia], [K]) of a vertex algebra $V$ 
is a linear endomorphism map $f$ of $V$ such that
\begin{eqnarray}
f(Y(u,x)v)=Y(f(u),x)v+Y(u,x)f(v)\;\;\;\mbox{ for }u,v\in V.
\end{eqnarray}
For any $u\in V$, $u_{0}$ is a derivation of $V$.
All derivations of $V$ form a Lie subalgebra $\Der V$ of $\gl(V)$
and all derivations $u_{0}$ for $u\in V$ form an ideal of $\Der V$,
which is a quotient algebra of the Lie algebra ${\cal{L}}(V)_{0}$.
For any Lie subalgebra $\g$ of $\Der V$, it is
straightforward to show (cf. [K]) that $V^{\g}$ is 
a vertex subalgebra of $V$.}
\er

Motivated by the Virasoro algebra and affine Lie algebras, 
following [DLM2] we next
consider certain quotient vertex algebras of ${\cal{V}}(R)$.
First we have (cf. [DLM2]):

\bl{lkeral-of-partial}
Let $R$ be a vertex Lie algebra and let $a\in R$ be such that
$\partial a=0$. Then
\begin{eqnarray}
& &Y_{-}(a,x)u=0=Y_{-}(u,x)a\;\;\;\mbox{ for all }u\in R\\
& &a(n)=0\;\;\;\mbox{ for }n\ne -1.
\end{eqnarray}
Furthermore, $a$ lies in the center of
${\cal{V}}(R)$ and $a(-1)$ is a central element of ${\cal{L}}(R)$.
\el

{\bf Proof.} With $\partial a=0$, we have
\begin{eqnarray}
{d\over dx}Y_{-}(a,x)=Y_{-}(\partial a,x)=0,
\end{eqnarray}
which immediately implies that $Y_{-}(a,x)=0$
because $Y_{-}(a,x)\in \Hom (R,x^{-1}R[x^{-1}])$.
It follows from the half skew symmetry
that $Y_{-}(u,x)a=0$ for $u\in R$. 

With $a_{i}b=0$ for $i\ge 0,\; b\in R$,
it follows from the commutator formula that $a(n)$ is 
a central element of ${\cal{L}}(R)$ for any $n\in \Z$.
We also have
\begin{eqnarray}
a(n)=\rho(a\otimes t^{n})={1\over n+1}\rho (\hat{\partial}(a\otimes
t^{n+1}))=0
\;\;\;\mbox{ for }n\ne -1,
\end{eqnarray}
as elements of ${\cal{L}}(R)$. Consequently,
\begin{eqnarray}
a(n)=0\;\;\;\mbox{ for }n\ne -1,
\end{eqnarray}
as operators on ${\cal{V}}(R)$. 
It follows from the commutator formula that $a$ is in the center of
${\cal{V}}(R)$.
$\;\;\;\;\Box$

For a vector space $U$, by {\em a partially defined linear functional} 
on $U$ we
mean a linear functional on some subspace of $U$. 
Rigorously speaking, a partially defined linear functional on $U$ 
consists of a linear subspace $D_{\lambda}$ of $U$ and 
a linear functional $\lambda$ on $D_{\lambda}$.

Let $\lambda$ be a partially defined linear functional on $\ker \partial$
with domain $D_{\lambda}$.
Denote by $I_{\lambda}$ the ${\cal{L}}(R)$-submodule of ${\cal{V}}(R)$
generated by the vectors $a-\lambda(a)$ 
for $a\in D_{\lambda}\subset \ker \partial$.
{}From Lemma \ref{lkeral-of-partial} we have
\begin{eqnarray}
{\cal{D}}(a-\lambda(a))={\cal{D}}a=a_{-2}{\bf 1}=a(-2){\bf 1}=0
\;\;\;\mbox{ for }a\in D_{\lambda}\subset \ker \partial,
\end{eqnarray}
where ${\cal{D}}$ is the ${\cal{D}}$-operator of ${\cal{V}}(R)$.
It follows that $I_{\lambda}$ is ${\cal{D}}$-stable.
Since $R$ generates ${\cal{V}}(R)$ as a vertex algebra,
$I_{\lambda}$ is an ideal of  ${\cal{V}}(R)$ (cf. [LL]). 
Following [DLM2] we set
\begin{eqnarray}
{\cal{V}}_{\lambda}(R)={\cal{V}}(R)/I_{\lambda},
\end{eqnarray}
a vertex algebra.
Note that ${\cal{V}}_{\lambda}(R)$ is still an ${\cal{L}}(R)$-module.
We also have
\begin{eqnarray}
{\cal{V}}_{\lambda}(R)
=U({\cal{L}}(R))\otimes_{U({\cal{L}}(R)_{+}\oplus D_{\lambda})}\C_{\lambda},
\end{eqnarray}
as an ${\cal{L}}(R)$-module, where $\C_{\lambda}=\C$ 
on which ${\cal{L}}(R)_{+}$ acts as zero
and $a$ acts as $\lambda(a)$ for $a\in D_{\lambda}\subset\ker\;\partial$,
where $D_{\lambda}$ is considered as a (central) subalgebra
of ${\cal{L}}(R)_{-}$ through the map $a\mapsto a(-1)$ 
(recall Lemma \ref{lkeral-of-partial}).

\br{r-general-lambda}
{\em Let $\lambda$ be the partially defined linear functional on 
$\ker_{R}\partial$ with $D_{\lambda}=0$. Then we easily see that
${\cal{V}}_{\lambda}(R)={\cal{V}}(R)$.}
\er

\br{r-on-primc}
{\em Let $R$ be a vertex algebra, let $\lambda$ 
be a partially defined linear functional on $\ker_{R}\partial$ and
let ${\cal{V}}_{\lambda}(R)$ be
the vertex algebra associated with $R$ and $\lambda$. 
Then it follows immediately from Primc's result (Theorem \ref{tprimc}) 
that for any vertex algebra $V$ and for any vertex Lie algebra 
homomorphism $f$ from $R$ to $V$ viewed as a vertex Lie algebra,
such that
\begin{eqnarray}
f(a)=\lambda(a){\bf 1}\;\;\;\mbox{ for }
a\in D_{\lambda}\subset\ker_{R}\partial,
\end{eqnarray}
there exists a unique vertex algebra homomorphism $\bar{f}$ from
${\cal{V}}_{\lambda}(R)$ to $V$ such that
$\bar{f}\psi$ extends $f$, where $\psi$ is the quotient map from
${\cal{V}}(R)$ onto ${\cal{V}}_{\lambda}(R)$.}
\er

On the other hand, denote by $J_{\lambda}$ the ideal of the symmetric
algebra $S(R)$, generated by the vectors 
$a-\lambda(a)$ for $a\in D_{\lambda}\subset\ker \partial$. Set
\begin{eqnarray}
S_{\lambda}(R)=S(R)/J_{\lambda},
\end{eqnarray}
the quotient algebra.
Since from Lemma \ref{lkeral-of-partial}
$u(i)(a-\lambda(a))=0$ for $u\in R,\; i\ge 0,\; 
a\in D_{\lambda}\subset\ker \partial$, 
$J_{\lambda}$ is an ${\cal{L}}(R)_{+}$-submodule of
$S(R)$, so that $S_{\lambda}(R)$ is also an ${\cal{L}}(R)_{+}$-module.

\bp{pquotientcase1}
Let $R$ be a vertex Lie algebra and let $\lambda$ be a 
partially defined linear functional on $\ker_{R} \partial$. Then
the linear map $\omega$ defined in Proposition \ref{psymmetrizationmap}
reduces to an ${\cal{L}}(R)_{0}$-isomorphism $\omega_{\lambda}$ from
$S_{\lambda}(R)$ onto ${\cal{V}}_{\lambda}(R)$. Furthermore,
for any Lie subalgebra $\g$ of ${\cal{L}}(R)_{0}$, the map $\omega$
reduces to a linear isomorphism from $S_{\lambda}(R)^{\g}$ onto
${\cal{V}}_{\lambda}(R)^{\g}$.
\ep

{\bf Proof.} Let $a\in D_{\lambda}\subset\ker \partial,
\; u^{(1)},\dots, u^{(n)}\in R$.
Since $[a(-1), u^{(i)}(-1)]=0$, we have
\begin{eqnarray*}
\omega( au^{(1)}\cdots u^{(n)})&=&
{1\over n!}\sum_{\sigma\in S_{n}} 
u^{(\sigma(1))}(-1)\cdots u^{(\sigma(n))}(-1)a(-1){\bf 1}\\
&=&{1\over n!}\sum_{\sigma\in S_{n}} u^{(1)}(-1)\cdots u^{(\sigma(n))}(-1)a,
\end{eqnarray*}
so that
\begin{eqnarray}
\omega((a-\lambda(a))u^{(1)}\cdots u^{(n)})&=&
{1\over n!}\sum_{\sigma\in S_{n}} 
u^{(\sigma(1))}(-1)\cdots u^{(\sigma(n))}(-1)(a-\lambda(a))\nonumber\\
&=&\left({1\over n!}\sum_{\sigma\in S_{n}} 
u^{(\sigma(1))}(-1)\cdots u^{(\sigma(n))}(-1)\right)
(a-\lambda(a)).
\end{eqnarray}
It follows that $\omega (J_{\lambda})=I_{\lambda}$, so $\omega$
gives rise to an ${\cal{L}}(R)_{0}$-module isomorphism $\omega_{\lambda}$ from
$S_{\lambda}(R)$ onto ${\cal{V}}_{\lambda}(R)$. The second assertion
follows easily.
$\;\;\;\;\Box$

\br{rusedlater}
{\em The relationship between $S(R)$ and ${\cal{V}}(R)$ and between
$S_{\lambda}(R)$ onto ${\cal{V}}_{\lambda}(R)$ will be 
further studied in Section 4 in the context of vertex Poisson algebras.}
\er

\section{Vertex Poisson algebras}
In this section we recall the notion of vertex Poisson algebra and 
we give certain general construction theorems of vertex Poisson algebras.
We apply the construction theorems to show that the symmetric algebra
over a vertex Lie algebra has a natural vertex Poisson algebra
structure, which was known to Frenkel and Ben-Zvi.

Throughout this paper, by a {\em differential algebra} we mean 
a commutative associative algebra $A$ with
identity $1$ equipped with a derivation $\partial$.
We also often denote the differential algebra by $(A,\partial)$.
We say that a subset $U$ of $A$ generates $A$ as a differential
algebra if $\partial^{n}U$ for $n\ge 0$ generate $A$ as an algebra.

The following notion of vertex Poisson algebra is due to [FB]
(cf. [DLM2] and [EF]):

\bd{dvpladefinition}
{\em A {\em vertex Poisson algebra} is a differential algebra
$(A,\partial)$ equipped with a vertex Lie algebra 
structure $(Y_{-},\partial)$ (with the same operator $\partial$)
such that for $a,b,c\in A$,
\begin{eqnarray}\label{evpaderivationaxiom}
Y_{-}(a,x)(bc)=(Y_{-}(a,x)b)c+b(Y_{-}(a,x)c).
\end{eqnarray}}
\ed

In terms of components, (\ref{evpaderivationaxiom}) amounts to
\begin{eqnarray}
a_{n}(bc)=(a_{n}b)c+b(a_{n}c)\;\;\;\mbox{ for }a,b,c\in A,\; n\ge 0,
\end{eqnarray}
where $Y_{-}(a,x)=\sum_{n\ge 0}a_{n}x^{-n-1}$.
That is, $a_{n}$, for all $a\in A,\; n\ge 0$, are derivations of $A$. 
Thus
\begin{eqnarray}
Y_{-}(a,x)\in x^{-1}(\Der\; A)[[x^{-1}]]\;\;\;\mbox{ for }a\in A.
\end{eqnarray}
This implies that $Y_{-}(a,x)1=0$, and then $Y_{-}(1,x)a=0$ by
half skew symmetry. Thus
\begin{eqnarray}\label{ey+1x=0}
Y_{-}(1,x)=0.
\end{eqnarray}

A vertex algebra $V$ 
is said to be {\em commutative} if $[Y(u,x_{1}),Y(v,x_{2})]=0$ 
for all $u,v\in V$. 
{}From [FHL], $V$ is commutative if and only if
$u_{n}v=0$ for all $u,v\in V,\; n\ge 0$.

\br{rcommutativeva}
{\em It was known (see [B1]; cf. [Li1], Proposition 2.1.6, [FB]) that 
if $V$ is a commutative vertex
algebra, then $V$ is an honest commutative associative algebra with
the product defined by
\begin{eqnarray}
u\cdot v=u_{-1}v\;\;\;\mbox{ for }u,v\in V
\end{eqnarray}
and with ${\bf 1}$ as the identity element. Furthermore,
the ${\cal{D}}$-operator ${\cal{D}}$ of $V$ is a derivation and
\begin{eqnarray}
Y(u,x)v=\left(e^{x{\cal{D}}}u\right)v\;\;\;\mbox{ for }u,v\in V.
\end{eqnarray}
Conversely, for any differential algebra $(A,\partial)$,
$(A,Y,1)$ carries the structure
of a commutative vertex algebra where $Y$ is defined by
$Y(a,x)b=(e^{x\partial}a)b$ for $a,b\in A$ (see [B1]).
This gives rise to a canonical isomorphism between
the category of commutative vertex algebras and
the category of differential algebras. In view of this,
a vertex Poisson algebra structure on a vector space $A$
consists of a commutative vertex algebra structure and a vertex Lie
algebra structure with a compatibility condition (see [FB]).}
\er

A {\em $\Z$-graded vertex Poisson algebra} is 
a vertex Poisson algebra $A$ equipped with a $\Z$-grading
$A=\coprod_{n\in \Z}A_{(n)}$ such that $A$ as an algebra is $\Z$-graded
and such that for $a\in A_{(n)},\; n,r\in \Z,\; m\in \N$,
\begin{eqnarray}
& &\partial A_{(m)}\subset A_{(m+1)},\\
& &a_{m}A_{(r)}\subset A_{(r+n-m-1)}.
\end{eqnarray}

An {\em ideal} of a vertex Poisson algebra $A$ is an ideal $I$ 
of $A$ as an associative algebra such that
\begin{eqnarray}
& &\partial I\subset I,\\
& &a_{n}I\subset I\;\;\;\mbox{ for }a\in A,\; n\ge 0.
\end{eqnarray}
{}From the half skew symmetry we also have
\begin{eqnarray}
u_{n}A\subset I\;\;\;\mbox{ for }u\in I,\; n\ge 0.
\end{eqnarray}
Thus, the quotient space $A/I$ has a natural
vertex Poisson algebra structure.

A {\em vertex-Poisson-algebra homomorphism} from $A$ to $B$
is an algebra homomorphism $f$ such that
\begin{eqnarray}
& &f\partial= f\partial\\
& &f(Y_{-}(a,x)b)=Y_{-}(f(a),x)f(b)\;\;\;\mbox{ for }a,b\in A.
\end{eqnarray}
A {\em vertex-Poisson-algebra isomorphism} is a bijective
vertex-Poisson-algebra homomorphism.

We have the following simple fact:

\bl{lsimplefacthomorphism}
Let $A$ and $B$ be vertex Poisson algebras and let $f$ be
an algebra homomorphism from $A$ to $B$ such that
$f\partial=\partial f$. Suppose that 
\begin{eqnarray}
f(Y_{-}(u,x)v)=Y_{-}(f(u),x)f(v)\;\;\;\mbox{ for }u,v\in U,
\end{eqnarray}
where $U$ is  a generating subset of $A$ as a differential algebra.
Then $f$ is a vertex-Poisson-algebra homomorphism.
\el

{\bf Proof.} Let $a,b\in A$ be such that
\begin{eqnarray}
f(Y_{-}(a,x)b)=Y_{-}(f(a),x)f(b).
\end{eqnarray}
Using the first equality of
(\ref{ethreetogether}) for both $A$ and $B$ and using
the assumption $f\partial=\partial f$ we obtain
\begin{eqnarray}
f(Y_{-}(a,x)\partial b)&=&f(\partial Y_{-}(a,x)b)-{d\over dx}f(Y_{-}(a,x)b)
\nonumber\\
&=&\partial f(Y_{-}(a,x)b)-{d\over dx}Y_{-}(f(a),x)f(b)
\nonumber\\
&=&\partial Y_{-}(f(a),x)f(b)-{d\over dx}Y_{-}(f(a),x)f(b)
\nonumber\\
&=&Y_{-}(f(a),x)\partial f(b)\nonumber\\
&=&Y_{-}(f(a),x)f(\partial b).
\end{eqnarray}
It follows from (\ref{evpaderivationaxiom}) and induction
that $$f(Y_{-}(u,x)a)=Y_{-}(f(u),x)f(a)
\;\;\;\mbox{ for any }u\in U,\; a\in A.$$
Using the half skew symmetry and the assumption $\partial f=f\partial$ 
we get
\begin{eqnarray}
f(Y_{-}(a,x)u)=Y_{-}(f(a),x)f(u)
\;\;\;\mbox{ for any }u\in U,\; a\in A.
\end{eqnarray}
Using (\ref{evpaderivationaxiom}) and induction again we get
$$f(Y_{-}(a,x)b)=Y_{-}(f(a),x)f(b)\;\;\;\mbox{ for any }a,b\in A.$$
Thus $f$ is a vertex Poisson algebra homomorphism.
$\;\;\;\;\Box$

In the following we prove certain general
results which will be useful for establishing
vertex Poisson algebra structures.

\bl{linductioncheckskewsymmetry}
Let $A$ be a differential algebra equipped with a linear map
$Y_{-}$ from $A$ to $x^{-1}(\Der\; A)[[x^{-1}]]$ such that 
$Y_{-}(u,x)v\in x^{-1}A[x^{-1}]$ for $u,v\in A$. 
Let $a,b,c,d\in A$ be such that
the half skew symmetry holds for all the pairs
$(a,c)$, $(a,d)$, $(b,c)$, $(b,d)$, $(ab,c)$, $(ab,d)$, $(a,cd)$, 
and $(b,cd)$. Then the half skew symmetry holds for the pair $(ab,cd)$:
\begin{eqnarray}\label{eabskewcd}
Y_{-}(ab,x)cd=\Sing\left( e^{x\partial}Y_{-}(cd,-x)ab\right).
\end{eqnarray}
\el

{\bf Proof.} Using all the assumptions we have
\begin{eqnarray}
& &Y_{-}(ab,x)(cd)\nonumber\\
&=&(Y_{-}(ab,x)c)d+cY_{-}(ab,x)d\nonumber\\
&=&\Sing \left( \left(e^{x\partial}Y_{-}(c,-x)(ab)\right)d+
ce^{x\partial}Y_{-}(d,-x)(ab)\right)\nonumber\\
&=&\Sing \left( (e^{x\partial}(Y_{-}(c,-x)a)b)d
+(e^{x\partial}(aY_{-}(c,-x)b))d\right)\nonumber\\
& &+\Sing \left( ce^{x\partial}(Y_{-}(d,-x)a)b
+ce^{x\partial}(aY_{-}(d,-x)b)\right)\label{emiddle}\\
&=&\Sing \left( (Y_{-}(a,x)c)(e^{x\partial}b)d+(e^{x\partial}a)(Y_{-}(b,x)c)d
\right)\nonumber\\
& &+\Sing \left(c(Y_{-}(a,x)d)(e^{x\partial}b)
+c(e^{x\partial}a)(Y_{-}(b,x)d)\right).\label{eyabcd}
\end{eqnarray}
Using (\ref{emiddle}) and the symmetry $(a,b,x)\rightarrow (c,d,-x)$
we get
\begin{eqnarray}
Y_{-}(cd,-x)(ab)
&=&\Sing \left( (e^{-x\partial}(Y_{-}(a,x)c)d)b
+(e^{-x\partial}(cY_{-}(a,x)d))b\right)\nonumber\\
& &+\Sing \left( ae^{-x\partial}(Y_{-}(b,x)c)d
+ae^{-x\partial}(cY_{-}(b,x)d)\right).
\end{eqnarray}
Then 
\begin{eqnarray}\label{eycdab}
& &\Sing \left(e^{x\partial}Y_{-}(cd,-x)(ab)\right)\nonumber\\
&=&\Sing \left( Y_{-}(a,x)c)d(e^{x\partial}b)
+(cY_{-}(a,x)d)(e^{x\partial}b)\right)\nonumber\\
& &+\Sing \left( (e^{x\partial}a)(Y_{-}(b,x)c)d
+(e^{x\partial}a)(cY_{-}(b,x)d)\right).
\end{eqnarray}
Combining (\ref{eyabcd}) with (\ref{eycdab}) we get (\ref{eabskewcd}).
$\;\;\;\;\Box$

Furthermore, we have:

\bp{pinductioncheckskewsymmetry}
Let $A$ be a differential algebra equipped with a linear map
$Y_{-}$ from $A$ to $x^{-1}(\Der\; A)[[x^{-1}]]$ such that 
$Y_{-}(u,x)v\in x^{-1}A[x^{-1}]$ for $u,v\in A$ and such that $Y_{-}(1,x)=0$. 
Suppose that
\begin{eqnarray}
Y_{-}(a,x)b=\Sing\left( e^{x\partial}Y_{-}(b,-x)a\right)
\end{eqnarray}
for $a\in A,\;b\in B$, where $B$ is a generating subset of $A$ 
as an algebra. 
Then the half skew symmetry holds for all
$a,b\in A$.
\ep

{\bf Proof.} Since $Y_{-}(a,x)\in x^{-1}(\Der\; A)[[x^{-1}]]$ for $a\in A$,
we have $Y_{-}(a,x)1=0$. Then
$Y_{-}(a,x)1=\Sing (e^{x\partial}Y_{-}(1,-x)a)\;(=0)$ for any $a\in A$.
Because $B$ generates $A$ as an algebra,
now it suffices to prove that 
\begin{eqnarray}
Y_{-}(a^{(1)}\cdots a^{(m)},x)(b^{(1)}\cdots b^{(n)})
=\Sing \left(e^{x\partial}
Y_{-}(b^{(1)}\cdots b^{(n)},-x)(a^{(1)}\cdots a^{(m)})\right)
\end{eqnarray}
for any $a^{(1)},\dots,a^{(m)}, b^{(1)},\dots, b^{(n)}\in B$ with
$m,n\ge 1$.
This follows immediately from Lemma \ref{linductioncheckskewsymmetry} 
(and induction on $m+n$).
$\;\;\;\;\Box$

Now we have the following basic result:

\bt{tinductioncheckPoisson}
Let $A$ be a differential algebra equipped with a linear map
$Y_{-}$ from $A$ to $x^{-1}(\Der\; A)[[x^{-1}]]$ such that 
$Y_{-}(a,x)b\in x^{-1}A[x^{-1}]$ for $a,b\in A$, $Y_{-}(1,x)=0$ and
\begin{eqnarray}\label{epartialderivativeintheorem}
Y_{-}(\partial a,x)={d\over dx}Y_{-}(a,x)\;\;\;\mbox{ for }a\in A.
\end{eqnarray}
Assume that $B$ is an ordered generating subset of $A$ 
as an algebra such that
for $a\in A,\; b\in B$,
\begin{eqnarray}
Y_{-}(a,x)b=\Sing\left( e^{x\partial}Y_{-}(b,-x)a\right)
\end{eqnarray}
and such that for $u,v,w\in B$ with $u\le v\le w$, 
\begin{eqnarray}\label{ehalfcommutatorinductionstep1}
[Y_{-}(u,x_{1}),Y_{-}(v,x_{2})]w
=\Sing \left( Y_{-}(Y_{-}(u, x_{1}-x_{2})v,x_{2})w\right).
\end{eqnarray}
Then $A$ is a vertex Poisson algebra.
The same assertion holds if we only assume that $B$ generates $A$ 
as a differential algebra and if in addition we assume 
\begin{eqnarray}\label{epartialbracketderivativeintheorem}
[\partial, Y_{-}(a,x)]={d\over dx}Y_{-}(a,x)\;\;\;\mbox{ for }a\in A.
\end{eqnarray}
\et

{\bf Proof.} By Proposition \ref{pinductioncheckskewsymmetry} we have
\begin{eqnarray}\label{eskewsymmetryinproof}
Y_{-}(a,x)b=\Sing\left( e^{x\partial}Y_{-}(b,-x)a\right)
\;\;\;\mbox{ for all }a,b\in A.
\end{eqnarray}
Furthermore, using (\ref{eskewsymmetryinproof}) and 
(\ref{epartialderivativeintheorem}) we get
\begin{eqnarray}
{d\over dx}Y_{-}(a,x)b&=&\Sing\left(\partial e^{x\partial}Y_{-}(b,-x)a
+e^{x\partial}{d\over dx}Y_{-}(b,-x)a\right)\nonumber\\
&=&\Sing\left(\partial e^{x\partial}Y_{-}(b,-x)a
-e^{x\partial}Y_{-}(\partial b,-x)a\right)\nonumber\\
&=&\partial Y_{-}(a,x)b-Y_{-}(a,x)\partial b.
\end{eqnarray}
Thus
\begin{eqnarray}\label{ebracket-derivative-in-proof}
[\partial, Y_{-}(a,x)]={d\over dx}Y_{-}(a,x)\;\;\;\mbox{ for any }
a\in A.
\end{eqnarray}
With the property (\ref{ebracket-derivative-in-proof}),
in view of Proposition \ref{ps3symmetry-commutator},
the half commutator formula for $(u,v,w)$ implies 
the half  commutator formula for any permutation of $(u,v,w)$.
Since $Y_{-}(a,x)\in x^{-1}(\Der\; A)[[x^{-1}]]$ and $B$ generates $A$,
(\ref{ehalfcommutatorinductionstep1}) holds for $u,v\in B$ and for
any $w\in A$. 
Again, since $Y_{-}(a,x)\in x^{-1}(\Der\; A)[[x^{-1}]]$ and $B$ generates $A$,
(\ref{ehalfcommutatorinductionstep1}) holds 
for $u\in B$ and for any $v, w\in A$.
Using this argument again we see that 
(\ref{ehalfcommutatorinductionstep1}) holds 
for any $u, v, w\in A$. This proves that $A$ is a vertex Poisson algebra.

For the second assertion, for $a\in A,\; b\in B$, 
using the conjugation formula (\ref{econjugation-formula}) (which
follows from (\ref{epartialderivativeintheorem}) and 
(\ref{epartialbracketderivativeintheorem})) we get
\begin{eqnarray}
Y_{-}(e^{x_{1}\partial}b,x)a
&=&e^{x_{1}{d\over dx}}Y_{-}(b,x)a\nonumber\\
&=&e^{x_{1}{d\over dx}}\Sing \left( e^{x\partial}
Y_{-}(a,-x)b\right)\nonumber\\
&=&\Sing_{x} \left( e^{(x+x_{1})\partial}e^{x_{1}{d\over dx}}
Y_{-}(a,-x)b\right)\nonumber\\
&=&\Sing_{x} \left( e^{x\partial}Y_{-}(a,-x)e^{x_{1}\partial}b\right).
\end{eqnarray}
Using the conjugation formula,
(\ref{ehalfcommutatorinductionstep1}) and the Taylor theorem 
we get
\begin{eqnarray}
& &[Y_{-}(e^{z_{1}\partial}u,x_{1}),Y_{-}(e^{z_{2}\partial}v,x_{2})]
e^{z\partial}w\nonumber\\
&=&e^{z_{1}{\partial\over\partial x_{1}}}
e^{z_{2}{\partial\over\partial x_{2}}}
[Y_{-}(u,x_{1}),Y_{-}(v,x_{2})]e^{z\partial}w\nonumber\\
&=&e^{(z_{1}-z){\partial\over\partial x_{1}}}
e^{(z_{2}-z){\partial\over\partial x_{2}}}
e^{z\partial}[Y_{-}(u,x_{1}),Y_{-}(v,x_{2})]w\nonumber\\
&=&e^{(z_{1}-z){\partial\over\partial x_{1}}}
e^{(z_{2}-z){\partial\over\partial x_{2}}}
e^{z\partial}\Sing\left( Y_{-}(Y_{-}(u,x_{1}-x_{2})v,x_{2})w\right)\nonumber\\
&=&e^{(z_{1}-z){\partial\over\partial x_{1}}}
e^{(z_{2}-z){\partial\over\partial x_{2}}}
\Sing_{x_{1},x_{2}}\left( Y_{-}(Y_{-}(u,x_{1}-x_{2})v,x_{2}+z)
e^{z\partial}w\right)\nonumber\\
&=&\Sing_{x_{1},x_{2}}
\left(Y_{-}(Y_{-}(u,x_{1}-x_{2}+z_{1}-z_{2})v,x_{2}+z_{2})
e^{z\partial}w\right)\nonumber\\ 
&=&\Sing_{x_{1},x_{2}}
\left(Y_{-}(Y_{-}(e^{(z_{1}-z_{2})\partial}u,x_{1}-x_{2})v,x_{2}+z_{2})
e^{z\partial}w\right)\nonumber\\ 
&=&\Sing_{x_{1},x_{2}}
\left(Y_{-}(e^{-z_{2}\partial}Y_{-}(e^{z_{1}\partial}u,x_{1}-x_{2})
e^{z_{2}\partial}v,x_{2}+z_{2})
e^{z\partial}w\right)\nonumber\\ 
&=&\Sing_{x_{1},x_{2}}
\left(Y_{-}(Y_{-}(e^{z_{1}\partial}u,x_{1}-x_{2})
e^{z_{2}\partial}v,x_{2})
e^{z\partial}w\right).
\end{eqnarray}
Since $\partial^{n}U$ for $n\ge 0$ generate $A$ as an algebra, 
now it follows immediately from the first assertion.
$\;\;\;\;\;\Box$

It is well known that for any Lie algebra $\g$, there exists a
unique Poisson algebra structure $\{\cdot,\cdot\}$ 
on the symmetric algebra $S(\g)$
such that $\{ u,v\}=[u,v]$ for $u,v\in \g$. 
On the other hand, if there is a Poisson algebra structure
$\{\cdot,\cdot\}$ on the symmetric algebra $S(U)$ for a vector space
$U$ such that $\{ u,v\}\in U$ for $u,v\in U$,  then 
$(U,\{\cdot,\cdot\})$ must be a Lie algebra. The following is 
a vertex analogue of this fact (cf. [FB], Example 15.2.2):

\bp{panother} Let $R$ be a vector space 
equipped with a linear operator $\partial$ and 
let $Y_{-}^{0}$ be a linear map from $R$ to $\Hom (R,x^{-1}R[x^{-1}])$.
Denote by $S(R)$ the symmetric algebra over $R$ and we extend
$\partial$ (uniquely) to a derivation of $S(R)$.
Then $Y_{-}^{0}$ extends to a vertex Poisson algebra structure $Y_{-}$ on
$(S(R),\partial)$ if and only if
$(R,\partial, Y_{-}^{0})$ carries the structure of a vertex Lie algebra.
Furthermore, such an extension is unique. 
\ep

{\bf Proof.} The ``only if'' part is clear. It is also clear that
if $Y_{-}^{0}$ extends to a vertex Poisson algebra structure $Y_{-}$ on
$(S(R),\partial)$, it must be unique. So, we must prove that
if $(R,\partial, Y_{-}^{0})$ carries the structure of a vertex Lie algebra,
then $Y_{-}^{0}$ extends to a vertex Poisson algebra structure $Y_{-}$
on $S(R)$.
First, for $u\in R$, we define a unique element 
\begin{eqnarray}
\tilde{Y}_{-}^{0}(u,x)\in x^{-1}(\Der \; S(R))[[x^{-1}]]
\end{eqnarray}
by
\begin{eqnarray}
\tilde{Y}_{-}^{0}(u,x)v=Y_{-}^{0}(u,x)v\;\;\;\mbox{ for }u,v\in R.
\end{eqnarray}
For $a\in S(R)$, we define
\begin{eqnarray}
Y_{-}(a,x)\in x^{-1}(\Der \; S(R))[[x^{-1}]]
\end{eqnarray}
by
\begin{eqnarray}
Y_{-}(a,x)u=\Sing
\left(e^{x\partial}\tilde{Y}_{-}^{0}(u,-x)a\right)\;\;\;\mbox{ for }u\in R.
\end{eqnarray}
For $u,v\in R$, using the half-skew symmetry of the vertex Lie algebra 
$(R,\partial,Y_{-}^{0})$ we have
\begin{eqnarray}
Y_{-}(v,x)u=\Sing\left(e^{x\partial}\tilde{Y}_{-}^{0}(u,-x)v\right)=
\Sing\left(e^{x\partial}Y_{-}^{0}(u,-x)v\right)=Y_{-}^{0}(v,x)u.
\end{eqnarray}
That is, $Y_{-}$ as a linear map from $R\otimes R$  
to $x^{-1}R[x^{-1}]$ extends $Y_{-}^{0}$. 
Since $Y_{-}(v,x),\; \tilde{Y}_{-}^{0}(v,x)\in x^{-1}(\Der\; S(R))[[x^{-1}]]$ and
$Y_{-}(v,x)u=\tilde{Y}_{-}^{0}(v,x)u$ for all $u\in R$, we have
\begin{eqnarray}
Y_{-}(v,x)a=\tilde{Y}_{-}^{0}(v,x)a\;\;\;\mbox{ for all }a\in S(R).
\end{eqnarray}
Similarly, since $[\partial,Y_{-}^{0}(u,x)]v={d\over dx}Y_{-}^{0}(u,x)v$ for $u,v\in R$,
it follows that for $u\in R,\; a\in S(R)$,
\begin{eqnarray}
[\partial,\tilde{Y}_{-}^{0}(u,x)]a={d\over dx}\tilde{Y}_{-}^{0}(u,x)a.
\end{eqnarray}
Then for $a\in S(R),\; u\in R$, we have
\begin{eqnarray}
Y_{-}(\partial a,x)u
&=&\Sing\left(e^{x\partial} \tilde{Y}_{-}^{0}(u,-x)\partial a\right)\nonumber\\
&=&\Sing\left(\partial e^{x\partial} \tilde{Y}_{-}^{0}(u,-x)a
+e^{x\partial}{d\over dx} \tilde{Y}_{-}^{0}(\partial u,-x)a
\right)\nonumber\\
&=&{d\over dx}\Sing\left( e^{x\partial} \tilde{Y}_{-}^{0}(u,-x)a\right)
\nonumber\\
&=&{d\over dx}Y_{-}(a,x)u.
\end{eqnarray}
Thus
$$Y_{-}(\partial a,x)={d\over dx}Y_{-}(a,x).$$

For $a\in S(R),\; u\in R$ we also have
\begin{eqnarray}
Y_{-}(a,x)u=\Sing \left( e^{x\partial}\tilde{Y}_{-}^{0}(u,-x)a\right)
=\Sing \left( e^{x\partial}Y_{-}(u,-x)a\right).
\end{eqnarray}
With $\tilde{Y}_{-}^{0}(u,x)v=Y_{-}^{0}(u,x)v$ for $u,v\in R$,
it now follows immediately from Theorem \ref{tinductioncheckPoisson}
that $(S(R),Y_{-},\partial)$ carries the structure of a vertex Poisson
algebra.
$\;\;\;\;\Box$

Furthermore, we have:

\bp{pwalgebra}
Let $R$ be a vertex Lie algebra and let
$\lambda$ be a partially defined linear functional on $\ker_{R} \partial$.
Then $J_{\lambda}$, the ideal of $S(R)$ generated by the
vectors $a-\lambda(a)$ for $a\in D_{\lambda}\subset\ker_{R} \partial$, 
is an ideal of
the vertex Poisson algebra $S(R)$, so that
the quotient space $S(R)/J_{\lambda}$, denoted by $S_{\lambda}(R)$,
is naturally a vertex Poisson algebra.
\ep

{\bf Proof.} 
For $a\in D_{\lambda}\subset\ker_{R} \partial$, 
we have $\partial (a-\lambda(a))=0$, so that
$\partial J_{\lambda}\subset J_{\lambda}$, since $\partial$ is a
derivation of $S(R)$. Similarly, 
\begin{eqnarray}
Y_{-}(u,x)J_{\lambda}\subset x^{-1}J_{\lambda}[[x^{-1}]]
\;\;\;\mbox{ for }u\in S(R),
\end{eqnarray}
since $Y_{-}(u,x)\in x^{-1}(\Der\; S(R))[[x^{-1}]]$ and 
$Y_{-}(u,x)(a-\lambda(a))=0$ for $a\in \ker_{R}\partial$.
(Recall Lemma \ref{lkeral-of-partial} that 
$Y_{-}(u,x)a=0$ for $a\in \ker_{R} \partial$.)
Therefore, $J_{\lambda}$ is an ideal of the vertex Poisson algebra $S(R)$.
$\;\;\;\;\Box$

Just like the Poisson algebra $S(\g)$ associated with a Lie algebra $\g$,
the vertex Poisson algebra $S(R)$ satisfies a certain
universal property.

\bp{puniversal-pva}
Let $R$ be a vertex Lie algebra and let $S(R)$ be the associated 
vertex Poisson algebra. Then for any vertex Poisson algebra $A$ and
for any vertex Lie algebra homomorphism $g$ from $R$ to $A$ viewed as
a vertex Lie algebra, there exists a unique vertex Poisson algebra
homomorphism $\bar{g}$ from $S(R)$ to $A$, extending $g$. 
Furthermore, if
\begin{eqnarray}
g(a)=\lambda(a)\;\;\;\mbox{ for }a\in D_{\lambda}\subset \ker_{R}\partial,
\end{eqnarray}
for some partially defined linear functional $\lambda$ on $\ker_{R}\partial$, 
the map $\bar{g}$ induces a vertex Poisson algebra
homomorphism from $S_{\lambda}(R)$ to $A$.
\ep

{\bf Proof.} The uniqueness is clear, since $R$ generates $S(R)$ as
a commutative associative algebra.
Let $\bar{g}$ be the unique algebra homomorphism from $S(R)$ to $A$ 
extending $g$. It remains to show 
\begin{eqnarray}
& &\bar{g} \partial (u)=\partial \bar{g} (u)\label{eg-partial-comm}\\
& &\bar{g}(u_{n}v)=\bar{g}(u)_{n}\bar{g}(v)\;\;\;\mbox{ for }u,v\in
S(R),\; n\ge 0.\label{eg-un-comm}
\end{eqnarray}
Let $u,v\in S(R)$ be such that $\bar{g} \partial (u)=\partial \bar{g}(u)$ and
$\bar{g}\partial (v)=\partial \bar{g}(v)$. Then
\begin{eqnarray}
\bar{g}(\partial (uv))
&=&\bar{g}(\partial (u))\bar{g}(v)+\bar{g}(u)\bar{g}(\partial (v))\nonumber\\
&=&(\partial \bar{g}(u))\bar{g}(v)+\bar{g}(u)\partial \bar{g}(v)\nonumber\\
&=&\partial \bar{g}(uv).
\end{eqnarray}
We also have
$\bar{g}\partial (1)=0=\partial \bar{g}(1)$ and 
\begin{eqnarray}
\bar{g}\partial (a)=g\partial a=\partial g(a)=\partial \bar{g}(a)
\;\;\;\mbox{ for }a\in R.
\end{eqnarray}
Thus $\ker (\bar{g}\partial-\partial \bar{g})$ is a subalgebra of $S(R)$,
containing $1$ and $R$, so that
$\ker (\bar{g}\partial-\partial \bar{g})=S(R)$. 
This proves (\ref{eg-partial-comm}).

Let $a\in R,\; n\in \N$. We have that $a_{n}$ is a derivation of $S(R)$,
$\bar{g}(a)_{n}$ is a derivation of $A$ and that
\begin{eqnarray}
& &\bar{g}(a_{n}1)=0=\bar{g}(a)_{n}\bar{g}(1)\\
& &\bar{g}(a_{n}b)=g(a_{n}b)=g(a)_{n}g(b)=\bar{g}(a)_{n}\bar{g}(b)
\;\;\;\mbox{ for }b\in R.
\end{eqnarray}
The same argument in the proof of (\ref{eg-partial-comm}) with $a_{n}$ in place of $\partial$
shows that
\begin{eqnarray}
\bar{g}(a_{n}u)=\bar{g}(a)_{n}\bar{g}(u)\;\;\;\mbox{ for }u\in S(R).
\end{eqnarray}
It follows from half skew symmetry and (\ref{eg-partial-comm}) that
\begin{eqnarray}
\bar{g}(u_{n}a)=\bar{g}(u)_{n}\bar{g}(a)
\;\;\;\mbox{ for }u\in S(R),\; n\in \N,\; a\in R.
\end{eqnarray}
Then (\ref{eg-un-comm}) follows from the same argument of the 
first paragraph. Therefore, $\bar{g}$ is a vertex Poisson algebra homomorphism.
The last assertion is clear.
$\;\;\;\;\Box$

Just as with the usual Poisson algebra $S(\g)$,
we call the vertex Poisson algebra structure on $S(R)$
associated with a vertex Lie algebra $R$
a {\em linear vertex Poisson structure}.
Next, we shall study ``nonlinear'' vertex Poisson structures on
a free differential algebra.

Let $U$ be a vector space space, fixed for the rest of this section.
Set
\begin{eqnarray}
A=S(\C[\partial]\otimes U),
\end{eqnarray}
the symmetric algebra over the space
$\C[\partial]\otimes U$, which is a free $\C[\partial]$-module.
We shall use the notation $\partial^{i}u$ for $\partial^{i}\otimes u$.
Extend $\partial$ (uniquely) to a derivation of $A$.
Now, $(A,\partial)$ is a differential algebra.
We refer to this differential algebra as the {\em free differential algebra 
over $U$}. 

We hope to determine all the vertex Poisson algebra structures on $A$.
Clearly, any vertex Poisson algebra structure $Y_{-}$ on $A$ is 
uniquely determined by giving $Y_{-}(u,x)v$ for $u,v\in U$.
Furthermore, the half skew-symmetry necessarily holds for any pair
$(u,v)\in U\times U$.

We define a {\em weak pre-vertex-Poisson structure} on $A$ 
to be a linear map $Y_{-}^{0}$ from $U\times U$ to $x^{-1}A[x^{-1}]$ 
such that 
\begin{eqnarray}
Y_{-}^{0}(u,x)v=\Sing\left(e^{x\partial }Y_{-}^{0}(v,-x)u\right)
\;\;\;\mbox{ for }u,v\in U.
\end{eqnarray}
The following result, to a certain extent, 
is analogous to Lemma 7.1 of [P]:

\bp{pweakprevla2}
Let $Y_{-}^{0}: (u,v)\mapsto Y_{-}^{0}(u,x)v$ be a weak pre-vertex-Poisson 
structure on $A$.
Then $Y_{-}^{0}$ extends uniquely to a linear map 
\begin{eqnarray}
Y_{-}: A&\rightarrow& \Hom (A,x^{-1}A[x^{-1}])\nonumber\\
       a&\mapsto& Y_{-}(a,x)=\sum_{n\in \N}a_{n}x^{-n-1}
\end{eqnarray}
such that for $a,b\in A$, 
\begin{eqnarray}
& &Y_{-}(a,x)=\sum_{n\in \N}a_{n}x^{-n-1}\in x^{-1}(\Der \;A)[[x^{-1}]]
\label{ederivationseries}\\
& &[\partial,Y_{-}(a,x)]=Y_{-}(\partial a,x)={d\over dx}Y_{-}(a,x)
\label{ebracket-partial-derivative}\\
& &Y_{-}(a,x)b=\Sing\left( e^{x\partial}Y_{-}(b,-x)a\right).
\label{ehalfskewsymmetryinstep}
\end{eqnarray}
\ep

{\bf Proof.} Suppose that $Y_{-}$ is such a linear map.
For $u,v\in U$, from (\ref{ebracket-partial-derivative}) we have
\begin{eqnarray}
Y_{-}(u,x)e^{x_{1}\partial}v=e^{x_{1}\partial}e^{-x_{1}{d\over dx}}Y_{-}(u,x)v
=e^{x_{1}\partial}e^{-x_{1}{d\over dx}}Y_{-}^{0}(u,x)v.
\end{eqnarray}
Since $\partial^{m}v$, for $m\in \N,\; v\in U$, generates $A$ as an
algebra and since $Y_{-}(u,x)\in x^{-1}(\Der\;A)[[x^{-1}]]$,
$Y_{-}(u,x)a$ for any $a\in A$ is uniquely determined.
Furthermore, using (\ref{ehalfskewsymmetryinstep})
and using the same argument we see that such linear map $Y_{-}$ 
is uniquely determined by the properties
(\ref{ederivationseries})-(\ref{ehalfskewsymmetryinstep}).

For the existence of a linear map $Y_{-}$ with the required
properties, first for $u\in U$ we define 
\begin{eqnarray}
\tilde{Y}_{-}^{0}(u,x)\in x^{-1}(\Der \;A)[[x^{-1}]]
\end{eqnarray}
in terms of generating functions by
\begin{eqnarray}\label{edefinitionfirststepgenerating}
\tilde{Y}_{-}^{0}(u,x)e^{x_{1}\partial}v
=e^{x_{1}\partial}e^{-x_{1}{d\over dx}}Y_{-}^{0}(u,x)v.
\end{eqnarray}
Clearly, $\tilde{Y}_{-}^{0}(u,x)\partial^{m}v\in
x^{-1}A[x^{-1}]$, so it follows immediately that 
$\tilde{Y}_{-}^{0}(u,x)\in \Hom (A,x^{-1}A[x^{-1}])$. It is also clear that
$\tilde{Y}_{-}^{0}(u,x)v=Y_{-}^{0}(u,x)v$ for $u,v\in U$.

Differentiating both sides of 
(\ref{edefinitionfirststepgenerating}) with respect to $x_{1}$ 
and using (\ref{edefinitionfirststepgenerating}) we get
\begin{eqnarray}
\tilde{Y}_{-}^{0}(u,x)\partial e^{x_{1}\partial}v
&=&\partial e^{x_{1}\partial}e^{-x_{1}{d\over dx}}Y_{-}^{0}(u,x)v
-e^{x_{1}\partial}{d\over dx}e^{-x_{1}{d\over dx}}Y_{-}^{0}(u,x)v\nonumber\\
&=&\partial \tilde{Y}_{-}^{0}(u,x)e^{x_{1}\partial}v
-{d\over dx} \tilde{Y}_{-}^{0}(u,x)e^{x_{1}\partial}v.
\end{eqnarray}
That is,
\begin{eqnarray}
[\partial,\tilde{Y}_{-}^{0}(u,x)]e^{x_{1}\partial}v
={d\over dx}\tilde{Y}_{-}^{0}(u,x)e^{x_{1}\partial}v.
\end{eqnarray}
Because $\partial,\;\tilde{Y}_{-}^{0}(u,x)\in x^{-1}(\Der\; A)[[x^{-1}]]$, we have
\begin{eqnarray}\label{ebracketderivativetildeB}
[\partial,\tilde{B}(u,x)]={d\over dx}\tilde{B}(u,x)
\;\;\;\mbox{for }u\in U.
\end{eqnarray}

Now, for $a\in A$, we define
\begin{eqnarray}
Y_{-}(a,x)=\sum_{n\in \N}a_{n}x^{-n-1}\in x^{-1}(\Der\; A)[[x^{-1}]]
\end{eqnarray}
by
\begin{eqnarray}	
Y_{-}(a,x)\partial^{m}u
=\Sing \left(e^{x\partial}\left(-{d\over dx}\right)^{m}
\tilde{Y}_{-}^{0}(u,-x)a\right)
\end{eqnarray}
for $m\in \N,\; u\in U$.
Clearly, $Y_{-}(a,x)\partial^{m}u\in x^{-1}A[x^{-1}]$.
Then
\begin{eqnarray}
Y_{-}(a,x)b\in x^{-1}A[x^{-1}]\;\;\;\mbox{ for }a,b\in A.
\end{eqnarray}

For $u,v\in U$, we have
\begin{eqnarray}
Y_{-}(u,x)v=\Sing\left( e^{x\partial}\tilde{Y}_{-}^{0}(v,-x)u\right)=
\Sing\left( e^{x\partial}Y_{-}^{0}(v,-x)u\right)=Y_{-}^{0}(u,x)v.
\end{eqnarray}
Hence $Y_{-}$ extends $Y_{-}^{0}$.

For $a\in A,\; u\in U,\; m\in \N$, we have
\begin{eqnarray}
& &[\partial,Y_{-}(a,x)] \partial^{m}u\nonumber\\
&=&\partial Y_{-}(a,x)\partial^{m}u-Y_{-}(a,x)\partial^{m+1}u\nonumber\\
&=&\Sing \left(
\partial e^{x\partial}\left(-{d\over dx}\right)^{m}\tilde{Y}_{-}^{0}(u,-x)a
-e^{x\partial}\left(-{d\over dx}\right)^{m+1}\tilde{Y}_{-}^{0}(u,-x)a\right)
\nonumber\\
&=&\Sing \left({d\over dx}\left\{e^{x\partial}\left(-{d\over dx}\right)^{m}
\tilde{Y}_{-}^{0}(u,-x)a\right\}\right)\nonumber\\
&=&{d\over dx}\Sing \left(\left\{e^{x\partial}\left(-{d\over dx}\right)^{m}
\tilde{Y}_{-}^{0}(u,-x)a\right\}\right)\nonumber\\
&=&{d\over dx}Y_{-}(a,x)\partial^{m}u,
\end{eqnarray}
and using (\ref{ebracketderivativetildeB}) we get
\begin{eqnarray}
Y_{-}(\partial a,x)\partial^{m}u
&=&\Sing \left(e^{x\partial}\left(-{d\over dx}\right)^{m}\tilde{Y}_{-}^{0}(u,-x)
\partial a\right)\nonumber\\
&=&\Sing \left(e^{x\partial}\left(-{d\over dx}\right)^{m}
\left(\partial +{d\over dx}\right)\tilde{Y}_{-}^{0}(u,-x)a\right)\nonumber\\
&=&{d\over dx}\Sing \left(e^{x\partial}\left(-{d\over dx}\right)^{m}
\tilde{Y}_{-}^{0}(u,-x)a\right)\nonumber\\
&=&{d\over dx}Y_{-}(a,x)\partial^{m}u.
\end{eqnarray}
Since 
$$\partial,\;\; Y_{-}(a,x),\;\; {d\over dx}Y_{-}(a,x),\;\;Y_{-}(\partial a,x)
\in x^{-1}(\Der\; A)[[x^{-1}]],$$
we immediately obtain (\ref{ebracket-partial-derivative}).
For $u,v\in U$, since $Y_{-}(u,x)v=\tilde{Y}_{-}^{0}(u,x)v$,
using (\ref{ebracket-partial-derivative}) we get
$$(Y_{-}(u,x)-\tilde{Y}_{-}^{0}(u,x))e^{x_{1}\partial}v
=e^{x_{1}\partial}(Y_{-}(u,x-x_{1})-\tilde{Y}_{-}^{0}(u,x-x_{1}))v=0.$$
Thus
\begin{eqnarray}
Y_{-}(u,x)=\tilde{Y}_{-}^{0}(u,x)\;\;\;\mbox{ on }A.
\end{eqnarray}

It now remains to prove the half skew symmetry. First, notice that 
using (\ref{ebracket-partial-derivative}) we get
\begin{eqnarray}
Y_{-}(a,x)\partial^{m}u
=\Sing \left(e^{x\partial}\left(-{d\over dx}\right)^{m}Y_{-}(u,-x)a\right)
=\Sing \left(e^{x\partial}Y_{-}(\partial^{m}u,-x)a\right).
\end{eqnarray}
Then it follows from Proposition \ref{pinductioncheckskewsymmetry} 
with $B=\C[\partial]U$ that the half skew symmetry holds for any pair
$(a,b)$ of elements of $A$.
$\;\;\;\;\Box$

The following theorem (which is somewhat
analogous to Proposition 7.7 of [P]) facilitates the construction of
vertex Poisson algebras:

\bt{tprevla}
Let $A$ be the free differential algebra over $U$ and let
$Y_{-}^{0}$ be a weak pre-vertex Poisson structure on $A$ such that
that for $u,v,w\in U$ from an ordered basis of $U$ with $u\le v\le w$,
\begin{eqnarray}\label{ebasicrelations}
& &\tilde{Y}_{-}^{0}(u,x_{1})Y_{-}^{0}(v,x_{2})w
-\tilde{Y}_{-}^{0}(v,x_{2})Y_{-}^{0}(u,x_{1})w\nonumber\\
&=&{\rm Sing}\left(e^{x_{2}\partial}
\tilde{Y}_{-}^{0}(w,-x_{2})Y_{-}^{0}(u,x_{1}-x_{2})v\right),
\end{eqnarray}
where $\tilde{Y}_{-}^{0}(u,x)\in x^{-1}(\Der\; A)[[x^{-1}]]$ is 
uniquely determined by 
\begin{eqnarray}
\tilde{Y}_{-}^{0}(u,x)e^{x_{1}\partial}v
=e^{x_{1}\partial}e^{-x_{1}{d\over dx}}Y_{-}^{0}(u,x)v
\end{eqnarray}
for $v\in U$.
Then $Y_{-}^{0}$ uniquely extends to a vertex Poisson structure $Y_{-}$ on 
$A$.
\et

{\bf Proof.} The uniqueness follows from Proposition \ref{pweakprevla2}.
Also, from Proposition \ref{pweakprevla2}, we have a linear map
$Y_{-}$ from $A$ to $\Hom (A,x^{-1}A[x^{-1}])$ such that
$Y_{-}(u,x)=\tilde{Y}_{-}^{0}(u,x)$ for $u\in U$ and such that
(\ref{ederivationseries})-(\ref{ehalfskewsymmetryinstep}) holds.
Then for $u,v,w\in U$, from (\ref{ebasicrelations}), using
(\ref{ehalfskewsymmetryinstep}) we get
\begin{eqnarray}
& &Y_{-}(u,x_{1})Y_{-}(v,x_{2})w-Y_{-}(v,x_{2})Y_{-}(u,x_{1})w\nonumber\\
&=&{\rm Sing}\left(e^{x_{2}\partial}
Y_{-}(w,-x_{2})Y_{-}(u,x_{1}-x_{2})v\right)\nonumber\\
&=&{\rm Sing}\left(Y_{-}(Y_{-}(u,x_{1}-x_{2})v,x_{2})w\right).
\end{eqnarray}
Now, it follows immediately from Theorem \ref{tinductioncheckPoisson}
that $Y_{-}$ is a vertex Poisson structure on $A$.
$\;\;\;\;\Box$

\bex{exaffinenonzeroc}
{\em Let $\g$ be a Lie algebra equipped with a nondegenerate symmetric 
invariant bilinear form $\<\cdot,\cdot\>$ and let $\ell$ be any
complex number.
Let $A$ be the free differential algebra over $\g$.
Define
\begin{eqnarray}
Y_{-}^{0}(u,x)v=[u,v] x^{-1}+\ell \<u,v\>x^{-2}\;\;\;\mbox{ for }u,v\in \g.
\end{eqnarray}
The bilinear map $B$ is a weak pre-vertex Poisson structure on $A$, since
$$Y_{-}^{0}(u,x)v=[v,u] (-x)^{-1}+\ell \<v,u\>(-x)^{-2}
=\Sing (e^{x\partial}B(v,-x)u).$$
Since $\tilde{B}(u,x)1=0$, we have
\begin{eqnarray}
& &\tilde{Y}_{-}^{0}(u,x_{1})Y_{-}^{0}(v,x_{2})w
-\tilde{Y}_{-}^{0}(v,x_{2})Y_{-}^{0}(u,x_{1})w\nonumber\\
&=&\tilde{Y}_{-}^{0}(u,x_{1})([v,w]x_{2}^{-1}+\ell \<v,w\>x_{2}^{-2})
-\tilde{Y}_{-}^{0}(v,x_{2})([u,w]x_{1}^{-1}+\ell \<u,w\> x_{1}^{-2})\nonumber\\
&=&Y_{-}^{0}(u,x_{1})[v,w]x_{2}^{-1}-Y_{-}^{0}(v,x_{2})[u,w]x_{1}^{-1}\nonumber\\
&=&[u,[v,w]]x_{1}^{-1}x_{2}^{-1}+\ell \<u, [v,w]\>x_{1}^{-2}x_{2}^{-1}
-[v,[u,w]]x_{1}^{-1}x_{2}^{-1}-\ell \<v, [u,w]\>x_{2}^{-2}x_{1}^{-1}
\nonumber\\
&=&[w, [u,v]]x_{1}^{-1}x_{2}^{-1}+\ell \<w,[u,v]\>
(x_{1}^{-1}x_{2}^{-2}+x_{1}^{-2}x_{2}^{-1}).
\end{eqnarray}
On the other hand,
\begin{eqnarray}
& &\Sing\left( e^{x_{2}\partial}\tilde{Y}_{-}^{0}(w,-x_{2})Y_{-}^{0}(u,x_{1}-x_{2})v\right)
\nonumber\\
&=& \Sing\left( e^{x_{2}\partial}
\left([w,[u,v]](-x_{2})^{-1}(x_{1}-x_{2})^{-1}
+\ell \<w,[u,v]\>(-x_{2})^{-2}(x_{1}-x_{2})^{-1}\right)\right)\nonumber\\
&=&[w,[u,v]] x_{1}^{-1}x_{2}^{-1}+\ell \<w,[u,v]\>
(x_{1}^{-1}x_{2}^{-2}+x_{1}^{-2}x_{2}^{-1}).
\end{eqnarray}
Thus
\begin{eqnarray}
\tilde{Y}_{-}^{0}(u,x_{1})Y_{-}^{0}(v,x_{2})w-\tilde{Y}_{-}^{0}(v,x_{2})Y_{-}^{0}(u,x_{1})w
=\Sing\left( e^{x_{2}\partial}\tilde{Y}_{-}^{0}(w,-x_{2})Y_{-}^{0}(u,x_{1}-x_{2})v\right)
\end{eqnarray}
for $u,v,w\in \g$. By Theorem \ref{tprevla}, 
$Y_{-}^{0}$ extends uniquely to a vertex Poisson structure on $A$.}
\eex

\br{raffinenonzeroc}
{\em From Example \ref{exaffinenonzeroc}, the subspace 
$R'=\C[\partial]\otimes \g \oplus \C$ of $A$ is a vertex Lie subalgebra
where $\partial \alpha =0$, $Y_{-}(u,x)\alpha=Y_{-}(\alpha,x)u=0$ 
and $Y_{-}(u,x)v=[u,v]x^{-1}+\ell \<u,v\>x^{-2}$ for $u,v\in \g$,
$\alpha\in \C$. It follows that there exists a unique
vertex Lie algebra structure on the space $R=\C[\partial]\otimes \g
\oplus \C c$, where $c$ is a symbol, 
such that $\partial c =0$, $Y_{-}(u,x)c=Y_{-}(c,x)u=0$ 
and $Y_{-}(u,x)v=[u,v]x^{-1}+c \<u,v\>x^{-2}$ for $u,v\in \g$.
This was proved in [P] (cf. [K]). }
\er

\section{Filtered vertex algebras and vertex Poisson algebras}

In this section we formulate a notion of filtered vertex algebra
and we prove that for any filtered vertex algebra $(V,E)$,
the associated graded vector space $\gr_{E}V$ is naturally a
vertex Poisson algebra. For any vertex algebra $V$, 
we determine all the filtered vertex algebras $(V,E)$ and we 
associate a canonical filtered vertex algebra $(V,E)$
to any $\N$-graded vertex algebra $V=\coprod_{n\ge 0}V_{(n)}$
with $V_{(0)}=\C {\bf 1}$. We also formulate and study a notion of 
generating subspace with PBW spanning property of a vertex algebra $V$
and we give a connection between 
this notion and our construction of ``good filtration'' of $V$.

\bd{dfilteredva}
{\em A {\em good filtration} of a vertex algebra $V$ 
is an increasing filtration by subspaces
\begin{eqnarray}
\cdots \subset E_{-2}\subset E_{-1}\subset E_{0}\subset E_{1}
\subset E_{2} \subset \cdots,
\;\;\; \cup_{n\in \Z}^{\infty}E_{n}=V,
\end{eqnarray}
such that ${\bf 1}\in E_{0}$ and such that
for $u\in E_{m},\; v\in E_{n}$,
\begin{eqnarray}
& &u_{j}v\in E_{m+n}\;\;\;\mbox{ for }j< 0\label{enegativemode}\\
& &u_{j}v\in E_{m+n-1}\;(\subset E_{m+n})\;\;\;\mbox{ for }j\ge 0.
\label{enonnegativemode}
\end{eqnarray}
A {\em filtered vertex algebra} is a vertex algebra
$V$ equipped with a good filtration $E=\{E_{n}\}$.
We sometimes denote the filtered vertex algebra by $(V,E)$.}
\ed

As an immediate  consequence of definition we have
\begin{eqnarray}
{\cal{D}}E_{m}\subset E_{m}\;\;\;\mbox{ for }m\in \Z,
\end{eqnarray}
since ${\cal{D}}v=v_{-2}{\bf 1}\in E_{m}$ for $v\in E_{m}$ by 
(\ref{enegativemode}) (notice that ${\bf 1}\in E_{0}$ by assumption).

Let $(V,E)$ be a filtered vertex algebra. 
Set
\begin{eqnarray}\label{egradedspace}
\gr_{E} V=\coprod_{n\in \Z}^{\infty}\left(\gr_{E}V\right)_{n}
=\coprod_{n\in \Z}^{\infty}(E_{n}/E_{n-1}),
\end{eqnarray}
a $\Z$-graded vector space.
In view of (\ref{enegativemode}), we have a bilinear 
multiplication ``$\cdot$'' on $A$ given by
\begin{eqnarray}\label{edefinitionmultiplication}
& &E_{m}/E_{m-1} \times E_{n}/E_{n-1}\rightarrow E_{m+n}/E_{m+n-1}
\nonumber\\
& &(u+E_{m-1})\cdot (v+E_{n-1})=u_{-1}v+E_{m+n-1}.
\end{eqnarray}
This makes $\gr_{E} V$ a non-associative algebra. It is clear that
${\bf 1}\in E_{0}$ is the identity.
Since ${\cal{D}}$ preserves $E_{m}$ for all $m$, we have a
linear operator $\partial$ on $\gr_{E} V$ defined by
\begin{eqnarray}\label{edefderivation}
\partial:  E_{m}/E_{m-1}&\rightarrow & E_{m}/E_{m-1}\nonumber\\
(u+E_{m-1})&\mapsto& {\cal{D}}u+E_{m-1}.
\end{eqnarray}
Furthermore, in view of (\ref{enonnegativemode}),
for every $i\ge 0$, we have a well defined bilinear multiplication on
$\gr_{E} V$:
\begin{eqnarray}\label{edefinitionvertexmultiplication}
E_{m}/E_{m-1}\times E_{n}/E_{n-1}&\rightarrow& E_{m+n-1}/E_{m+n-2}
\nonumber\\
(u+E_{m-1},v+E_{n-1})&\mapsto& u_{i}v+E_{m+n-2}.
\end{eqnarray}
Then for $u\in E_{m},\; v\in E_{n}$ with $m,n\in \Z$ 
we define 
\begin{eqnarray}
Y_{-}(u+E_{m-1},x)(v+E_{n-1})=\sum_{i\ge 0}(u_{i}v+E_{m+n-2})x^{-i-1}
\in x^{-1}(E_{m+n-1}/E_{m+n-2})[x^{-1}].
\end{eqnarray}
We have:

\bp{pfilteredvoa}
Let $(V,E)$ be a filtered vertex algebra and let $\gr_{E}V$ be
the associated graded vector space defined in (\ref{egradedspace}).
Then $\gr_{E} V$ equipped with the structures $(\cdot, \partial, Y_{-})$
defined in (\ref{edefinitionmultiplication}), (\ref{edefderivation})
and (\ref{edefinitionvertexmultiplication})
is a vertex Poisson algebra. Furthermore, 
\begin{eqnarray}
& & \partial (\gr_{E}V)_{m}\subset (\gr_{E}V)_{m}\\
& &ab\in (\gr_{E}V)_{m+n}=E_{m+n}/E_{m+n-1}\\
& &a_{i}b\in (\gr_{E} V)_{m+n-1}=E_{m+n-1}/E_{m+n-2}
\end{eqnarray}
for $m,n\in \Z,\; a\in E_{m}/E_{m-1},\; b\in E_{n}/E_{n-1},\; i\ge 0$.
\ep

{\bf Proof.} First, we prove that $(\gr_{E}V,\cdot)$ is a
commutative associative algebra with ${\bf 1}+E_{-1}$ as the identity.
Clearly, ${\bf 1}+E_{0}$ is the identity for the algebra.
For $u\in E_{m},\; v\in E_{n},\; m,n\in \Z$,
using (\ref{ecommutatorincomponents}), (\ref{enegativemode}) and 
(\ref{enonnegativemode}), we have
\begin{eqnarray}
u_{-1}v-v_{-1}u=u_{-1}v_{-1}{\bf 1}-v_{-1}u_{-1}{\bf 1}=
\sum_{i\ge 0}{-1\choose i}(u_{i}v)_{-2-i}{\bf 1}\in E_{m+n-1}.
\end{eqnarray}
This implies that $\gr_{E} V$ is commutative.
Furthermore, let $u\in E_{m},\; v\in E_{n},\; w\in E_{k}$.
Using the iterate formula (\ref{eiterateincomponents}) we have
\begin{eqnarray}\label{especialiterate}
(u_{-1}v)_{-1}w=\sum_{i\ge 0}(-1)^{i}{-1\choose i} 
\left( u_{-1-i}v_{-1+i}w-(-1)^{-1}v_{-2-i}u_{i}w\right).
\end{eqnarray}
Notice that from (\ref{enonnegativemode}) and (\ref{enegativemode})
we have $u_{-1-i}v_{-1+i}w\in E_{m+n+k-1}$ for $i\ge 1$
and $v_{-2-i}u_{i}w\in E_{m+n+k-1}$ for $i\ge 0$.
Then from (\ref{especialiterate}) we get
\begin{eqnarray}
(u_{-1}v)_{-1}w+E_{m+n+k-1}=u_{-1}v_{-1}w+E_{m+n+k-1}.
\end{eqnarray}
This proves that $\gr_{E} V$ is associative, so that
$\gr_{E} V$ is a commutative associative algebra.

We next prove that the operator $\partial$ is a derivation.
Let $u\in E_{m},\; v\in E_{n}$ with $m,n\in \Z$. 
Using (\ref{edbracketformulaincomponenets}) we get
\begin{eqnarray}
& &\partial \left((u+E_{m-1})(v+E_{n-1})\right)\nonumber\\
&=&{\cal{D}}(u_{-1}v)+E_{m+n-1}\nonumber\\
&=&u_{-1}{\cal{D}}v+({\cal{D}}u)_{-1}v+E_{m+n-1}\nonumber\\
&=&(u+E_{m-1})\partial (v+E_{n-1})+(\partial (u+E_{m-1}))(v+E_{n-1}).
\end{eqnarray}
This proves that $\partial$ is a derivation of the commutative associative
algebra $\gr_{E} V$.

Now we check the axioms involving $Y_{-}$. 
Let $u\in E_{m},\; v\in E_{n},\;w\in E_{k}$. Then
\begin{eqnarray}
& &Y_{-}(u+E_{m-1},x)\left((v+E_{n-1})(w+E_{k-1})\right)\nonumber\\
&=&\sum_{i\ge 0}(u_{i}v_{-1}w+E_{m+n+k-2})x^{-i-1}\nonumber\\
&=&\sum_{i\ge 0}(v_{-1}u_{i}w+E_{m+n+k-2})x^{-i-1}
+\sum_{i,j\ge 0}{i\choose j}((u_{j}v)_{i-1-j}w
+E_{m+n+k-2})x^{-i-1}\nonumber\\
&=&\sum_{i\ge 0}(v_{-1}u_{i}w+(u_{i}v)_{-1}w+E_{m+n+k-2})x^{-i-1}
\nonumber\\
& &+\sum_{i\ge 0}\sum_{0\le j<i}{i\choose j}((u_{j}v)_{i-1-j}w
+E_{m+n+k-2})x^{-i-1}\nonumber\\
&=&\sum_{i\ge 0}(v_{-1}u_{i}w+(u_{i}v)_{-1}w+E_{m+n+k-2})x^{-i-1}
\nonumber\\
&=&(v+E_{n-1})Y_{-}(u+E_{m-1},x)(w+E_{k-1})
+\left(Y_{-}(u+E_{m-1},x)(v+E_{n-1})\right)(w+E_{k-1}),\nonumber
\end{eqnarray}
since for $0\le j<i$,
\begin{eqnarray}
(u_{j}v)_{i-1-j}w\in E_{m+n+k-2}.
\end{eqnarray}
This proves that $Y_{-}(u+E_{m-1},x)\in x^{-1}(\Der\; \gr_{E} V)[[x^{-1}]]$.
Thus
\begin{eqnarray}
Y_{-}(a,x)\in x^{-1}(\Der\; \gr_{E} V)[[x^{-1}]]\;\;\;\mbox{ for }a\in \gr_{E} V.
\end{eqnarray}
We also have
\begin{eqnarray}
& &\partial Y_{-}(u+E_{m-1},x)(v+E_{n-1})\nonumber\\
&=&\sum_{i\ge 0}({\cal{D}}u_{i}v+E_{m+n-2})x^{-i-1}\nonumber\\
&=&\sum_{i\ge 0}(u_{i}{\cal{D}}v+({\cal{D}}u)_{i}v+E_{m+n-2})x^{-i-1}
\nonumber\\
&=&Y_{-}(u+E_{m-1},x)\partial (v+E_{n-1})
+Y_{-}(\partial (u+E_{m-1}),x)(v+E_{n-1}).
\end{eqnarray}
Using the formula ${\cal{D}}u_{i}v=u_{i}{\cal{D}}v-iu_{i-1}v$ instead,
we get
\begin{eqnarray}
& &\partial Y_{-}(u+E_{m-1},x)(v+E_{n-1})\nonumber\\
&=&\sum_{i\ge 0}({\cal{D}}u_{i}v+E_{m+n-2})x^{-i-1}\nonumber\\
&=&\sum_{i\ge 0}(u_{i}{\cal{D}}v-iu_{i-1}v+E_{m+n-2})x^{-i-1}\nonumber\\
&=&Y_{-}(u+E_{m-1},x)\partial (v+E_{n-1})
+{d\over dx}Y_{-}(u+E_{m-1},x)(v+E_{n-1}).
\end{eqnarray}
This proves that 
$[\partial, Y_{-}(u+E_{m-1},x)]=Y_{-}(\partial (u+E_{m-1}),x)
={d\over dx}Y_{-}(u+E_{m-1},x)$.
Thus
\begin{eqnarray}
[\partial, Y_{-}(a,x)]=Y_{-}(\partial a,x)={d\over dx}Y_{-}(a,x)
\;\;\;\mbox{ for all }a\in \gr_{E} V.
\end{eqnarray}

Finally, for $u\in E_{m},\; v\in E_{n-1},\; w\in E_{k}$, using
(\ref{ecommutatorincomponents}) we have
\begin{eqnarray}
& &[Y_{-}(u+E_{m-1},x_{1}), Y_{-}(v+E_{n-1},x_{2})](w+E_{k-1})\nonumber\\
&=&\sum_{p,q\ge 0} \left(u_{p}v_{q}w-v_{q}u_{p}w+E_{m+n+k-3}\right)
x_{1}^{-p-1}x_{2}^{-q-1}\nonumber\\
&=&\sum_{p,q\ge 0}\sum_{i=0}^{p}{p\choose i}
\left( (u_{i}v)_{p+q-i}w+E_{m+n+k-3}\right)
x_{1}^{-p-1}x_{2}^{-q-1}\nonumber\\
&=&\sum_{p,q\ge 0}\sum_{i=0}^{p}{p\choose i}
\left( (u+E_{m-1})_{i}(v+E_{n-1})\right)_{p+q-i}(w+E_{k-1})
x_{1}^{-p-1}x_{2}^{-q-1},
\end{eqnarray}
where $Y_{-}(u+E_{m-1},x)=\sum_{p\ge 0}(u+E_{m-1})_{p}x^{-p-1}$
with $(u+E_{m-1})_{p}\in \Der\; \gr_{E}V$.
This proves the half commutator formula for a vertex Lie algebra.
Now the proof is complete. $\;\;\;\;\Box$

A good filtration $E=\{E_{n}\}$ of a vertex algebra $V$ 
is said to be {\em truncated} if $E_{n}=0$ for $n$ sufficiently small.
Typical truncated good filtrations are those with $E_{n}=0$ for $n<0$.
(Note that for any good filtration $E$ of a nonzero vertex algebra $V$,
$E_{0}\ne 0$, since ${\bf 1}\in E_{0}$.) 

\bp{psimple-vas}
Let $V$ be a simple vertex algebra in the sense that there is no 
ideal other than $V$ and $0$. Then for any truncated good filtration of $V$ 
we have $E_{n}=0$ for $n<0$.
\ep

{\bf Proof.} Suppose that $E=\{E_{n}\}$ is a good filtration of $V$ such that
$E_{-k}\ne 0$ and $E_{n}=0$ for $n<-k$, where $k$ is a positive integer.
For $u,v\in E_{-k}$, from (\ref{enegativemode}) and
(\ref{enonnegativemode}) we have
\begin{eqnarray}
u_{m}v\in E_{-2k}=0\;\;\;\mbox{ for }m\in \Z.
\end{eqnarray}
In particular,
\begin{eqnarray}
u_{m}u=0\;\;\;\mbox{ for }u\in E_{-k},\; m\in \Z.
\end{eqnarray}
Thus  each element $u$ of $E_{-k}$ is nilpotent in the sense of [LL].
By Proposition 3.10.6 of [LL], all nilpotent elements of $V$ form an
ideal, which must equal $V$ because $E_{-k}\ne 0$ and $V$ is simple. 
This is a contradiction because ${\bf 1}$ is not nilpotent.
$\;\;\;\;\Box$

\br{rfiltration-E0}
{\em Let $E=\{E_{n}\}$ be a good filtration of a vertex algebra $V$.
For $u,v\in E_{0}$, from (\ref{enegativemode}) and
(\ref{enonnegativemode}) we have
\begin{eqnarray}
& &u_{n}v\in E_{-1}\;(\subset E_{0})\;\;\;\mbox{ for }n\ge 0,\label{eideal-1}\\
& &u_{n}v\in E_{0}\;\;\;\mbox{ for }n< 0.
\end{eqnarray}
Consequently, $E_{0}$ is a vertex subalgebra of $V$ 
(notice that ${\bf 1}\in E_{0}$ by assumption).
Furthermore, $E_{-1}$ is an ideal of $E_{0}$, since
for $u\in E_{0},\; v\in E_{-1}$, 
\begin{eqnarray}
u_{n}v,\; v_{n}u\in E_{-1}\;\;\;\mbox{ for }n\in \Z.
\end{eqnarray}
It follows from (\ref{eideal-1}) that $E_{0}/E_{-1}$ is a commutative vertex algebra.
In particular, if $E_{-1}=0$, then $E_{0}$ is a commutative vertex subalgebra of $V$.}
\er

We say a filtration $E=\{E_{n}\}$ for $V$ is {\em finer} than
$E'=\{E_{n}'\}$, or $E'$ is {\em coarser} than $E$,
if $E_{n}\subset E_{n}'$ for all $n\in \Z$.

\br{rtrivialfiltration}
{\em For any vertex algebra $V$, we have a {\em trivial} filtration
$E=\{E_{n}\}$ with $E_{n}=0$ for $n<0$, $E_{0}=\C {\bf 1}$ and $E_{n}=V$ for $n\ge 1$.
It is clear that this is a good filtration of $V$.
Of course, it is the coarsest filtration. 
For this trivial filtration, we have 
$\gr_{E}V=E_{0}\oplus E_{1}/E_{0}=\C {\bf 1}\oplus V/\C {\bf 1}$,
where $(V/\C {\bf 1})\cdot (V/\C {\bf 1})=0$
and $\gr_{E}V=\C{\bf 1}\oplus V/\C {\bf 1}$ as a vertex Lie algebra
is the direct sum of
vertex Lie algebras $\C$ and $V/\C{\bf 1}$, the quotient vertex Lie
algebra of $V$ by the ideal $\C {\bf 1}$ (cf. [K]).}
\er

The following theorem gives a general construction and classification 
of good filtrations with $E_{n}=0$ for $n<0$ and $E_{0}=\C {\bf 1}$
for a vertex algebra:

\bt{tclassification}
Let $V$ be any vertex algebra. Let $U$ be a subspace of
$V$ equipped with a vector space decomposition $U=\coprod_{n\ge 1}U_{n}$.
For $n\ge 0$, define $E_{n}^{U}$ to be the subspace of $V$ 
linearly spanned by the vectors 
\begin{eqnarray}
u^{(1)}_{-n_{1}}\cdots u^{(r)}_{-n_{r}}{\bf 1}
\end{eqnarray}
for $r\ge 0,\; u^{(i)}\in U_{m_{i}},\; n_{i}\ge 1$ with $m_{i}\ge 1$ and
$m_{1}+\cdots +m_{r}\le n$.
Suppose that $U$ generates $V$ as a vertex algebra, i.e.,
\begin{eqnarray}\label{egeneratingspannweak}
V={\rm span}\{
u^{(1)}_{n_{1}}\cdots u^{(r)}_{n_{r}}{\bf 1}\;|\;
\mbox{ for }r\ge 0,\; u^{(i)}\in U,\; n_{i}\in \Z\},
\end{eqnarray}
and suppose that \begin{eqnarray}\label{ebasicinduction}
u_{i}v\in E^{U}_{m+n-1}\;\;\;\mbox{ for }u\in U_{m},\; v\in
U_{n},\;i\ge 0,\;m,n\ge 1,
\end{eqnarray}
i.e., for $u\in U_{m},\; v\in
U_{n},\;i\ge 0,\;m,n\ge 1$, 
$u_{i}v$ is a linear combination of the vectors 
$$u^{(1)}_{-n_{1}}\cdots u^{(r)}_{-n_{r}}{\bf 1}$$
for $r\ge 0,\; u^{(i)}\in U_{m_{i}},\; n_{i}\ge 1$ with
$m_{1}+\cdots +m_{r}\le m+n-1$.
Then $E_{U}=\{E_{n}^{U}\}$ is a good filtration of $V$. In particular,
\begin{eqnarray}\label{egeneratingspann1}
V=\cup_{n\ge 0}E_{n}^{U}=
{\rm span}\{
u^{(1)}_{-n_{1}}\cdots u^{(r)}_{-n_{r}}{\bf 1}\;|\;
\mbox{ for }r\ge 0,\; u^{(i)}\in U,\; n_{i}\ge 1\}.
\end{eqnarray}
Furthermore, any good filtration $E=\{E_{n}\}$ of $V$ with $E_{n}=0$ for $n<0$
and with $E_{0}=\C {\bf 1}$ can be obtained this way.
\et

{\bf Proof.} First, we prove that for any good filtration
$E=\{E_{n}\}$ with the special property, 
there exists a subspace $U$ of $V$, equipped with a vector space 
decomposition $U=\coprod_{n\ge 1}U_{n}$, such that $E=E_{U}$
and such that 
(\ref{egeneratingspann1}) and (\ref{ebasicinduction}) hold.
For $n\ge 1$, let $U_{n}$ be a subspace of $E_{n}$ such that
$E_{n}=U_{n}\oplus E_{n-1}$. Set $U=\coprod_{n\ge 1}U_{n}$.
(It is easy to see that the sum is direct.)
We are going to prove that $E_{n}=E_{n}^{U}$ for all $n\ge 0$.
Since $U_{m}\subset E_{m}$ for $m\ge 1$, it follows from
the definition of $E_{n}^{U}$ and (\ref{enonnegativemode}) 
that $E_{n}^{U}\subset E_{n}$ for $n\ge 1$.
By definition, we also have $E_{0}=\C {\bf 1}=E_{0}^{U}$.
{}From definition, we have $U_{n}\subset E_{n}^{U}$ 
(since $u=u_{-1}{\bf 1}$) and $E_{n}=U_{n}+E_{n-1}$
for $n\ge 1$. It follows from induction that $E_{n}\subset E_{n}^{U}$,
so $E=E_{U}$. 
Since $\cup_{n\ge 0}E_{n}^{U}=\cup_{n\ge 0}E_{n}=V$,
(\ref{egeneratingspann1}) holds. Furthermore, 
(\ref{ebasicinduction}) holds because
$u_{i}v\in E_{m+n-1}=E^{U}_{m+n-1}$.

Let $U$ be a subspace of $V$ equipped with a vector space 
decomposition $U=\prod_{n\ge 1}U_{n}$ such that 
(\ref{egeneratingspann1}) and (\ref{ebasicinduction}) hold.
{}From now on we shall just use $E_{n}$ for $E_{n}^{U}$.
To prove that $E=\{E_{n}=E_{n}^{U}\}$ is a good filtration 
we must prove that $\cup_{n\ge 0}E_{n}=V$ and that
for $n,p\ge 0,\; m\in \Z$ and for $v\in E_{p}$, 
\begin{eqnarray}
& &v_{m}E_{n}\subset E_{p+n-1}\;\;\;\;\;\;\;\;\mbox{ for }m\ge 0
\label{e3.18}\\
& &v_{m}E_{n}\subset E_{p+n}\;\;\;\;\;\;\;\;\;\;\;\mbox{ for }m< 0.
\label{e3.19}
\end{eqnarray}
Notice that (\ref{e3.18}) and (\ref{e3.19}) imply
that $\cup_{n\ge 0}E_{n}$ is a vertex subalgebra of $V$.
Since $U\subset \cup_{n\ge 0}E_{n}$ and $U$ generates $V$,
(\ref{e3.18}) and (\ref{e3.19}) imply that $\cup_{n\ge 0}E_{n}=V$.
Then it suffices to prove (\ref{e3.18}) and (\ref{e3.19}).

For $u\in U_{n}$ with $n\ge 1$, we say $u$ is {\em homogeneous} 
and we define $|u|=n$. First, we observe that from definition,
\begin{eqnarray}\label{efirstfact}
u_{m}E_{n}\subset E_{n+|u|}\;\;\;\mbox{ for homogeneous }u\in
U,\;\mbox{ and for } m<0.
\end{eqnarray}
We also have
\begin{eqnarray}
{\cal{D}}E_{n}\subset E_{n}\;\;\;\mbox{ for }n\ge 0,
\end{eqnarray}
since $[{\cal{D}},u_{-k}]=ku_{-k-1}$ and ${\cal{D}}{\bf 1}=0$ 
for $u\in U,\; k\ge 1$.

We shall prove (\ref{e3.18}) and (\ref{e3.19}) by induction on $n+p$.
If $p=0$ or $n=0$, it is clear because $E_{0}=\C {\bf 1}$
and $v_{-k-1}{\bf 1}={1\over k!}{\cal{D}}^{k}v\in E_{p}$ for $v\in
E_{p},\;k\ge 0$.
Assume that $n,p\ge 1$. We only need to consider
typical vectors 
$$v=v^{(1)}_{-k_{1}}\cdots v^{(r)}_{-k_{r}}{\bf 1}$$
for $r\ge 1,\; v^{(i)}\in U,\; k_{i}\ge 1$ with
$|v^{(1)}|+\cdots+|v^{(r)}|=p$.

If $r\ge 2$, we have
\begin{eqnarray*}
v=u_{-k}v'\;\;\;\mbox{ where }u=v^{(1)}\in U,\; 
k=k_{1}\ge 1,\; v'\in E_{p-|u|}
\end{eqnarray*}
with $|u|<p$.
For $w\in E_{n}$, by the iterate formula we have
\begin{eqnarray}
v_{m}w=\sum_{i\ge 0}{-k\choose i}(-1)^{i}
\left( u_{-k-i}v'_{m+i}w-(-1)^{k+i}v'_{-k+m-i}u_{i}w\right).
\end{eqnarray}
Noticing that $u\in E_{|u|}, \; v'\in E_{p-|u|}$ and $|u|+n<p+n$,
for $i\ge 0$, using the inductive hypothesis and (\ref{efirstfact}) we have
\begin{eqnarray}
& &u_{-k-i}v'_{m+i}w-(-1)^{k+i}v'_{-k+m-i}u_{i}w\in E_{n+p-1}
\;\;\;\mbox{ for }m\ge 0\\
& &u_{-k-i}v'_{m+i}w-(-1)^{k+i}v'_{-k+m-i}u_{i}w\in E_{n+p}
\;\;\;\mbox{ for }m< 0.
\end{eqnarray}
Then it follows that $v_{m}w\in E_{n+p-1}$ for $m\ge 0$
and $v_{m}w\in E_{n+p}$ for $m< 0$.

Now, we consider the case that $r=1$. 
That is, $v=u_{-k}{\bf 1}$ for some $u\in U,\;
k\ge 1$ with $p=|u|$, so
$$v_{m}=(u_{-k}{\bf 1})_{m}={m\choose k-1}u_{m+1-k}
\;\;\;\mbox{ for }m\in \Z.$$ 
If $m<0$, then $m+1-k<0$, and if $0\le m\le k-1$, then
$v_{m}=0$ and if $m\ge k$, we have $m+1-k\ge 0$.
In view of this, it suffices to consider $v=u$.
Consider any vector $u^{(1)}_{-n_{1}}w\in E_{n}$ with $w\in E_{n-|u^{(1)}|}$.
Then
\begin{eqnarray}
u_{m}u^{(1)}_{-n_{1}}w=u^{(1)}_{-n_{1}}u_{m}w
+\sum_{i\ge 0}{m\choose i}(u_{i}u^{(1)})_{m-n_{1}-i}w.
\end{eqnarray}
Since $u\in E_{|u|},\; w\in E_{n-|u^{(1)}|}$ and $|u|+n-
|u^{(1)}|=p+n-|u^{(1)}|<n+p$, 
by inductive hypothesis, $u_{m}w\in E_{|u|+n-|u^{(1)}|}$ for $m<0$
and $u_{m}w\in E_{|u|+n-|u^{(1)}|-1}$ for $m\ge 0$, so that
$u^{(1)}_{-n_{1}}u_{m}w\in E_{|u|+n}$ for $m<0$
and $u^{(1)}_{-n_{1}}u_{m}w\in E_{|u|+n-1}$ for $m\ge 0$.
On the other hand, since, by assumption
$$u_{i}u^{(1)}\in E_{|u|+|u^{(1)}|-1},$$
by inductive hypothesis, we have $(u_{i}u^{(1)})_{m-n_{1}-i}w\in
E_{|u|+n-1}$ for all $i\ge 0,\; m\in \Z$.
Thus
$u_{m}u^{(1)}_{-n_{1}}w\in E_{n+p}$ for $m<0$ and
$u_{m}u^{(1)}_{-n_{1}}w\in E_{n+p-1}$ for $m\ge 0$.
This completes the inductive step and the proof.$\;\;\;\;\Box$

For an illustration we next apply Theorem \ref{tclassification} to
the vertex algebra associated with a vertex Lie algebra.
Recall from Proposition \ref{panother}
that for any vertex Lie algebra $R$,
the symmetric algebra $S(R)$ has a vertex Poisson algebra structure. Furthermore,
$S(R)$ is an ${\cal{L}}(R)_{+}$-module.

\bp{pisomorphism}
Let $R$ be a vertex Lie algebra and let ${\cal{V}}(R)$ be the associated
vertex algebra. Let $E=\{E_{n}\}$ be the sequence defined by
\begin{eqnarray}
E_{n}={\rm span}\{ a^{(1)}(-1)\cdots a^{(r)}(-1){\bf
1}\;|\; 0\le r\le n,\; a^{(i)}\in R\}. 
\end{eqnarray}
Then $E=\{E_{n}\}$ is a good filtration of ${\cal{V}}(R)$ and each
$E_{n}$ is an ${\cal{L}}(R)_{+}$-submodule of ${\cal{V}}(R)$.
Furthermore, the linear map 
\begin{eqnarray}
\psi:S(R)&\rightarrow & \gr_{E} {\cal{V}}(R)\nonumber\\
u^{(1)}\cdots u^{(n)}&\mapsto& 
u^{(1)}(-1)\cdots u^{(n)}(-1){\bf 1}+E_{n-1}
\end{eqnarray}
for $u^{(i)}\in R$ is a vertex-Poisson-algebra isomorphism. Furthermore,
$\psi$ is an ${\cal{L}}(R)_{+}$-module isomorphism.
\ep

{\bf Proof.} Taking $U=U_{1}=R$ we have $E_{n}=E_{n}^{U}$ as defined in 
Theorem \ref{tclassification}. Since $R$ generates ${\cal{V}}(R)$ as a vertex algebra
and since $u_{i}v\in R$ for $u,v\in R,\;
i\ge 0$, by Theorem \ref{tclassification}
$E$ is a good filtration. By (\ref{e3.18}) and (\ref{e3.19}) we have
$$v_{m}E_{n}\subset E_{n}\;\;\;\;\mbox{ for }v\in R,\; m,n\ge 0.$$
Then each $E_{n}$ is an ${\cal{L}}(R)_{+}$-submodule of ${\cal{V}}(R)$.

It follows from Proposition \ref{pprimc} and
the Poincar\'e-Birkhoff-Witt theorem that
$\psi$ is a linear isomorphism (see [Di]) and then it follows from
the definition of the multiplication of $\gr_{E}{\cal{V}}(R)$ that
$\psi$ is an algebra isomorphism.
Furthermore, for $u,v\in R$, we have
\begin{eqnarray}
\psi (\partial u)=\partial u+E_{0}
={\cal{D}}u+E_{0}=\partial (u+E_{0})
=\partial \psi (u),
\end{eqnarray}
and
\begin{eqnarray}
\psi (Y_{-}(u,x)v)=\sum_{i\ge 0}(u_{i}v+E_{0})x^{-i-1}
=Y_{-}(u+E_{0},x)(v+E_{0})=Y_{-}(\psi (u),x)f(v).
\end{eqnarray}
By Lemma \ref{lsimplefacthomorphism},
$\psi$ is a vertex-Poisson-algebra homomorphism, so that it is a
vertex-Poisson-algebra isomorphism.

Note that the ${\cal{L}}(R)_{+}$-module structure on $S(R)$ is
given by $u(x)^{-}=Y_{-}(u,x)$ for $u\in R$, where
$$u(x)^{-}=\sum_{n\ge 0}u(n)x^{-n-1}\in {\cal{L}}(R)_{+}[[x^{-1}]]$$
and $Y_{-}$ is the vertex Lie algebra structure map of $S(R)$.
On the other hand, the ${\cal{L}}(R)_{+}$-module structure 
on ${\cal{V}}(R)$ is
given by $u(x)^{-}=Y(u,x)^{-}$ for $u\in R$, where $Y$ is the vertex
operator map of ${\cal{V}}(R)$ and $Y(u,x)^{-}=\Sing Y(u,x)$,
so that
the ${\cal{L}}(R)_{+}$-module structure on $\gr_{E}{\cal{V}}(R)$ is
given by $u(x)^{-}=Y_{-}(u+E_{0},x)$ for $u\in R=U\subset E_{1}$, 
where $Y_{-}$ is the
vertex Lie algebra structure map of $\gr_{E}{\cal{V}}(R)$.
Now, with $\psi$ as a vertex-Poisson-algebra isomorphism, 
it follows immediately that $\psi$ is 
an ${\cal{L}}(R)_{+}$-module isomorphism.
$\;\;\;\;\Box$

Let $U$ be a subspace of a vertex algebra $V$.
Following [K] we say that {\em $U$ strongly generates $V$} if
\begin{eqnarray}
V={\rm span}\{ u^{(1)}_{-n_{1}}\cdots u^{(r)}_{-n_{r}}{\bf
1}\;|\; r\ge 0,\; u^{(i)}\in U,\; n_{i}\ge 1\}.
\end{eqnarray}
We say that {\em $U$ generates $V$ with PBW spanning property}
if for some basis
$\{u^{(\alpha)}\;|\;\alpha\in I\}$ of $U$ and for any order ``$>$''
on the set $\{ u^{(\alpha)}_{-r}\;|\; \alpha\in I,\; r\ge 1\}$,
$V$ is linearly spanned by the vectors
\begin{eqnarray}
u^{(\alpha_{1})}_{-n_{1}}\cdots u^{(\alpha_{r})}_{-n_{r}}{\bf 1}
\end{eqnarray}
for $r\ge 0,\; \alpha_{i}\in I,\; n_{i}\ge 1$, with
$u^{(\alpha_{1})}_{-n_{1}}\ge \cdots \ge u^{(\alpha_{r})}_{-n_{r}}$.
In terms of this notion we have:

\bt{tgeneral-pbw-spanning-property}
Let $V$ be a vertex algebra and let $U$ be a generating subspace of $V$. 
Suppose that there exists a vector-space decomposition
$U=\coprod_{n\ge 1}U_{n}$ such that for any $u\in U_{m},\; v\in
U_{n},\; i\ge 0$, $u_{i}v$ is a linear combination of the vectors 
$$u^{(1)}_{-n_{1}}\cdots u^{(r)}_{-n_{r}}{\bf 1}$$
for $r\ge 0,\; u^{(i)}\in U_{m_{i}},\; n_{i}\ge 1$ with
$m_{1}+\cdots +m_{r}\le m+n-1$.
Then $U$ generates $V$ with PBW spanning property.
\et

{\bf Proof.} By Theorem \ref{tclassification}, we have a good filtration
$E=\{E_{n}^{U}\}$ associated to $U$. In particular,
$$V=\cup_{n\ge 0}E_{n}^{U}
={\rm span}\{ u^{(1)}_{-n_{1}}\cdots u^{(r)}_{-n_{r}}{\bf
1}\;|\; r\ge 0,\; u^{(i)}\in U,\; n_{i}\ge 1\}.$$
Let $\{ u^{(\alpha)}\;|\; \alpha\in I\}$ be a basis
of $U$, consisting of homogeneous vectors. For $\alpha,\beta\in I,\; r,s\in \Z$, from
Borcherds' commutator formula we have
\begin{eqnarray}
[u^{(\alpha)}_{r}, u^{(\beta)}_{s}]=\sum_{i\ge 0} {r\choose i}
(u^{(\alpha)}_{i}u^{(\beta)})_{r+s-i}.
\end{eqnarray}
Furthermore, if $u^{(\alpha)}\in U_{m}\;(\subset E_{m}^{U})$ and 
$u^{(\beta)}\in U_{n}\;(\subset E_{n}^{U})$, 
by (\ref{enonnegativemode}) we have
$$u^{(\alpha)}_{i}u^{(\beta)}\in E^{U}_{m+n-1}\;\;\;\mbox{ for }i\ge 0.$$
With the property (\ref{enegativemode}),
it now follows immediately from the classical argument
with the universal enveloping algebra of a Lie algebra (cf. [Di]), just as in [KL].
$\;\;\;\;\Box$

\br{rnonlinearly-generated-vas}
{\em Notice that the vertex algebra ${\cal{V}}(R)$ associated with
a vertex Lie algebra $R$ is generated by $U=R$ with {\em linear} relations
while Theorems \ref{tclassification} and \ref{tgeneral-pbw-spanning-property}
concern vertex algebras generated by $U$ with 
general {\em nonlinear} generating relations.}
\er

We next apply Theorems \ref{tclassification} and 
\ref{tgeneral-pbw-spanning-property} to ``$\N$-graded vertex algebras.''

A {\em $\Z$-graded vertex algebra} is a vertex algebra
equipped with a $\Z$-grading $V=\coprod_{n\in \Z}V_{(n)}$ such that
the following conditions hold for $u\in V_{(k)},\; m,n,k\in \Z$:
\begin{eqnarray}
u_{m}V_{(n)}\subset V_{(n+k-m-1)}.
\end{eqnarray}
An {\em $\N$-graded vertex algebra} is defined in the obvious way.
For $v\in V_{(n)}$ for $n\in \Z$, we say $v$ is {\em homogeneous of weight $n$}
and we write $\wt v=n$.

The following results were proved in [KL] (cf. [Li3]):

\bt{tkl1}
Let $V=\coprod_{n\in \N}V_{(n)}$ be an $\N$-graded vertex algebra
with $V_{(0)}=\C {\bf 1}$.
Denote by $C_{1}(V)$ the subspace of $V_{+}=\coprod_{n>0}V_{(n)}$ 
linearly spanned by the vectors $u_{-1}v$ for $u,v\in V_{+}$ 
and by the vectors $w_{-2}{\bf 1}\;(={\cal{D}}w)$ for $w\in V$. 
Then a graded subspace $U$ of $V_{+}$ strongly generates $V$, i.e.,
$$V={\rm span}\{
u^{(1)}_{-n_{1}}\cdots u^{(r)}_{-n_{r}}{\bf 1}
\;|\; r\ge 0,\; u^{(i)}\in U,\; n_{i}\ge 1\},$$ 
if and only if $V_{+}=U+C_{1}(V)$. Furthermore,
$U$ is a minimal graded strong generating subspace of $V$ if and only
if $V_{+}=U\oplus C_{1}(V)$.
If a graded subspace $U$ of $V_{+}$ strongly generates $V$, 
then $U$ generates $V$ with PBW spanning property.
\et

Now we have:

\bt{tequivalent-generating-property}
Let $V=\coprod_{n\in \N}V_{(n)}$ be an $\N$-graded vertex algebra
with $V_{(0)}=\C {\bf 1}$ and
let $U$ be a graded subspace of $V_{+}=\coprod_{n>0}V_{(n)}$. Then
the following statements are equivalent:

(a) $U$ generates $V$ as a vertex algebra and 
for $u\in U_{m}=U\cap V_{(m)},\; v\in U_{n}=U\cap V_{(n)},\; m,n\ge 1$ 
and for $i\ge 0$,
$u_{i}v$ is a linear combination of the vectors 
$$u^{(1)}_{-n_{1}}\cdots u^{(r)}_{-n_{r}}{\bf 1}$$
for $r\ge 0,\; u^{(i)}\in U_{m_{i}},\; n_{i}\ge 1$ with
$m_{1}+\cdots +m_{r}\le m+n-1$.

(b) $U$ strongly generates $V$, i.e.,
\begin{eqnarray}\label{egeneralvectors}
V={\rm span}\{ u^{(1)}_{-k_{1}}\cdots u^{(r)}_{-k_{r}}{\bf 1}\;|\;
r\ge 0,\; u^{(i)}\in U,\; k_{i}\ge 1\}.
\end{eqnarray}

(c) $V=U+C_{1}(V)$.

(d) $U$ generates $V$ with PBW spanning property.

Furthermore, assume that any one of the four equivalent conditions holds and
for $n\ge 0$, denote by $E_{n}^{U}$ the subspace of $V$ linearly spanned
by the vectors 
$$u^{(1)}_{-k_{1}}\cdots u^{(r)}_{-k_{r}}{\bf 1}$$
for $r\ge 0,\; u^{(i)}\in U,\; k_{i}\ge 1$
with $\wt u^{(1)}+\cdots +\wt u^{(r)}\le n$. 
Then the sequence $E_{U}=\{E_{n}^{U}\}$ is a good filtration
of vertex algebra $V$ and for $n\ge 0$,
\begin{eqnarray}\label{esumsubset}
V_{(0)}+\cdots+V_{(n)}\subset E_{n}^{U}.
\end{eqnarray}
\et

{\bf Proof.} From Theorem \ref{tkl1}, (b) and (c) are equivalent, and (b) implies
(d). Clearly, (d) implies (b). Thus (b), (c) and (d) are equivalent.
In view of Theorem \ref{tclassification} (note that $U=\coprod_{n\ge 1}(U\cap V_{(n)})$), 
(a) implies (b). Also in view of Theorem \ref{tclassification},
it suffices to prove that (b) implies (a) and (\ref{esumsubset}).
First we prove that (b) implies (\ref{esumsubset}).
{}From our assumption, for $m\ge 0$,
 $V_{(m)}$ is linearly spanned by the vectors
$$u^{(1)}_{-k_{1}}\cdots u^{(r)}_{-k_{r}}{\bf 1}$$
for homogeneous vectors $u^{(i)}\in U$ with
$$m=(\wt u^{(1)}+k_{1}-1)+\cdots +(\wt u^{(r)}+k_{r}-1).$$
Since $k_{i}\ge 1$, we have $m\ge \wt u^{(1)}+\cdots +\wt u^{(r)}$, 
so that by definition,
$$u^{(1)}_{-k_{1}}\cdots u^{(r)}_{-k_{r}}{\bf 1}\in E_{m}^{U}.$$
Thus $V_{(m)}\subset E_{m}^{U}$. We immediately have (\ref{esumsubset}), 
since $E_{m}^{U}\subset E_{n}^{U}$ for $m\le n$.

For any $u\in U_{(m)},\; v\in U_{(n)},\; m,n\ge 1,\;i \ge 0$,
since $\wt (u_{i}v)=\wt u+\wt v-i-1=m+n-i-1$ and $V_{(m)}\subset
E_{m}^{U}$ for all $m\ge 0$, we have
$$u_{i}v\in V_{(m+n-i-1)}\subset E_{m+n-i-1}^{U}\subset E_{m+n-1}^{U}.$$
This proves that (b) implies (a), completing the proof. $\;\;\;\;\Box$

\br{rnotunique}
{\em  Suppose that $E=\{E_{n}\}$ is a good filtration of $V$. We define
$E'=\{ E'_{n}\}$ by $E'_{2n}=E_{n}$ and $E_{2n+1}'=E_{n}$ 
for $n\ge 0$. Clearly, $E'$ is an increasing filtration of $V$.
It is straightforward to check that it is a good filtration.
Thus, in general there does not exist a finest good filtration.
Of course, good filtrations of a vertex algebra $V$ by no means
are unique.}
\er

Despite of the non-uniqueness of good filtrations on a vertex algebra
$V$ we are going to show that the good filtration $E^{U}=\{E_{n}^{U}\}$ 
associated with a graded subspace $U$ in fact does not depend on $U$. 
First, we have:

\bl{lgoodfiltration}
Let $V=\coprod_{n\ge 0}V_{(n)}$ be an $\N$-graded vertex algebra
with $V_{(0)}=\C {\bf 1}$ and let 
$E=\{E_{n}\}$ be the filtration obtained in 
Theorem \ref{tequivalent-generating-property} from
a graded subspace $U$ of $V_{+}$. 
Let $E'=\{E'_{n}\}$ be any good filtration
for $V$ such that $V_{(n)}\subset E'_{n}$ for $n\ge 0$.
Then 
\begin{eqnarray}
E_{n}\subset E'_{n}\;\;\;\mbox{ for all }n\ge 0.
\end{eqnarray}
\el

{\bf Proof.} By definition, $E_{0}=\C {\bf 1}\subset E'_{0}$.
Assume $n\ge 1$. By definition, $E_{n}$ is linearly spanned by the vectors
$$u^{(1)}_{-m_{1}}\cdots u^{(r)}_{-m_{r}}{\bf 1}$$
for $r\ge 1,\; u^{(i)}\in U,\; m_{i}\ge 1$ with
$\wt u^{(1)}+\cdots + \wt u^{(r)}\le n.$
For $1\le i\le r$, by assumption 
$$u^{(i)}\in V_{(\wt u^{(i)})}\subset E'_{\wt u^{(i)}}.$$
In view of (\ref{enegativemode}) we have
$$u^{(1)}_{-n_{1}}\cdots u^{(r)}_{-n_{r}}{\bf 1}\in E'_{\wt
u^{(1)}+\cdots +\wt u^{(r)}}\subset E'_{n}.$$
Thus $E_{n}\subset E'_{n}$. $\;\;\;\;\Box$.

Since for the filtration $E^{U}=\{E_{n}^{U}\}$ associated to 
any graded strong generating subspace $U$, $V_{(n)}\subset E_{n}$
for all $n\ge 0$ (by Theorem \ref{tequivalent-generating-property}), 
it follows from Lemma \ref{lgoodfiltration} that the filtrations
associated to any two graded generating subspaces must be the same.
Therefore, we have proved:

\bt{tklfilteredvoa}
Let $V=\coprod_{n\in \N}V_{(n)}$ be an $\N$-graded vertex algebra 
with $V_{(0)}=\C {\bf 1}$.
Then the filtration $E_{U}=\{E_{n}^{U}\}$ constructed in
Theorem \ref{tequivalent-generating-property} from
a graded strong generating subspace $U$ of $V$ does not depend on $U$.
Furthermore, this filtration is the unique finest filtration with
the property that $V_{(n)}\subset E_{n}$ for all $n\ge 0$.
$\;\;\;\;\Box$
\et

We call the good filtration associated with the graded subspace 
$U=V_{+}$ the {\em standard filtration} and
we call the associated vertex Poisson algebra
the {\em standard vertex Poisson algebra associated with $V$},
which we denote by $\gr\; V$. In a sequel we shall use 
$\gr\; V$ to study the vertex algebraic structure of $V$.

\section{Formal deformation of vertex (Poisson) algebras}
In this section we first formulate a notion of $h$-adic vertex algebra
and then using this notion we formulate 
a notion of formal deformation of vertex algebras and vertex Poisson algebras.
We relate the construction of vertex Poisson algebras from 
filtered vertex algebras with Frenkel and Ben-Zvi's construction
and we give a formal deformation of the vertex Poisson algebras 
$S(R)$ and $S_{\lambda}(R)$ associated with a vertex Lie algebra $R$.

In the literature, vertex algebras are often considered to be
over $\C$ while most of the results naturally carry over for 
vertex algebras over a field of characteristic $0$, or over 
a unital commutative associative algebra over
a field of characteristic $0$.
In this section we shall concern vertex algebras
over $\C$, or over a unital commutative associative algebra $K$ 
over $\C$, e.g., $K=\C [h]$, the polynomial algebra for
a formal variable $h$.

Let $V$ be a vertex algebra over $\C$. Then $V[h]=V\otimes \C[h]$
is naturally a vertex algebra over $\C[h]$ and 
$V\otimes \C[[h]]$ is naturally a vertex algebra over $\C[[h]]$.
Consider the space $V[[h]]$ of all formal power series in $h$.
(Notice that $V\otimes \C[[h]]$ is a proper subspace of $V[[h]]$
unless $V$ is finite-dimensional.)
It is well known that $V[[h]]$ is the completion of 
the $\C[h]$-module $V[h]$ with respect to the $h$-adic topology. 
Formally extend the vertex operator map $Y$ of vertex algebra $V$ 
to a $\C[[h]]$-linear map from $V[[h]]$ to 
$(\End\; V[[h]])[[x,x^{-1}]]$ by
\begin{eqnarray}
Y(a(h),x)b(h)=\sum_{m,n\ge 0}Y(a(m),x)b(n) h^{m+n}\;
\left(\in V[[h]][[x,x^{-1}]]\right)
\end{eqnarray}
for $a(h)=\sum_{n\ge 0}a(n)h^{n},\; b(h)=\sum_{n\ge 0}b(n)h^{n}\in V[[h]]$.
Because $Y(a(h),x)b(h)$ in general involves 
infinitely many negative powers of $x$,
the triple $(V[[h]],Y, {\bf 1})$ does not carry the structure 
of a vertex algebra (over $\C[h]$ or $\C[[h]]$) in the precise sense.
On the other hand, note that for any nonnegative integer $n$,
the quotient space $V[[h]]/h^{n}V[[h]]$ is naturally a vertex algebra
(over $\C[[h]]$, or $\C[h]$), which is isomorphic to 
the quotient vertex algebra $V[h]/h^{n}V[h]$ over $\C[h]$.
Motivated by this we formulate the following notion of 
$h$-adic vertex algebra:

\bd{d-h-adicvoa}
{\em An {\em $h$-adic vertex algebra} is a $\C[[h]]$-module
$A[[h]]$, where $A$ is a vector space over $\C$,
equipped with the $h$-adic topology,
a distinguished vector ${\bf 1}\in A[[h]]$ 
and a continuous $\C[[h]]$-linear map $Y_{h}$ from 
$A[[h]]\otimes A[[h]]$ to $A[[h]][[x,x^{-1}]]$
such that for every nonnegative integer $n$,
the triple $(A[[h]]/h^{n}A[[h]],\bar{Y}_{(n)},{\bf 1}+h^{n}A[[h]])$
carries the structure of a vertex algebra over $\C[[h]]$,
where $\bar{Y}_{(n)}$ is the natural quotient map of $Y_{h}$, i.e., 
\begin{eqnarray}
\bar{Y}_{(n)}(a+h^{n}A[[h]],x)(b+h^{n}A[[h]])
=\sum_{m\in \Z}(a_{m}b+h^{n}A[[h]])x^{-m-1}
\end{eqnarray}
for $a,b\in A$, where $Y_{h}(a,x)b=\sum_{m\in \Z}a_{m}b x^{-m-1}$.}
\ed

\br{rinverse-limits}
{\em Let $A[[h]]$ be an $h$-adic vertex algebra.
For every pair of nonnegative integers
$i$ and $j$ with $i\le j$, it is clear that
the natural map $\theta_{ij}$ from $A[[h]]/h^{j}A[[h]]$ 
onto $A[[h]]/h^{i}A[[h]]$ is a
vertex algebra homomorphism over $\C[[h]]$. 
The vertex algebras $A[[h]]/h^{n}A[[h]]$ (over $\C[[h]]$)
together with these vertex algebra morphisms form
an inverse system of vertex algebras over $\C[[h]]$.
Then the $h$-adic vertex algebra $A[[h]]$ can be thought of as
an inverse limit of the vertex algebras $A[[h]]/h^{n}A[[h]]$ 
over $\C[[h]]$.}
\er

\bp{pexplicit-version}
Let $A$ be a vector space over $\C$. An $h$-adic vertex algebra
structure on $A[[h]]$ amounts to a distinguished vector ${\bf 1}\in A[[h]]$ 
and a continuous $\C[[h]]$-linear map $Y_{h}$ from 
$A[[h]]\otimes A[[h]]$ to $A[[h]][[x,x^{-1}]]$
such that the following conditions hold: For $a\in A[[h]]$,
\begin{eqnarray}
& &Y_{h}({\bf 1},x)a=a,\\
& &Y_{h}(a,x){\bf 1}\in A[[h]][[x]]\;\;\mbox{ and }
\;\;\lim_{x\rightarrow 0}Y_{h}(a,x){\bf 1}=a;
\end{eqnarray}
for $a,b\in A[[h]]$, $n\in \N$, 
there exists an integer $k$ such that
\begin{eqnarray}\label{e-hadic-truncation}
a_{m}b\in h^{n}A[[h]]\;\;\;\mbox{ for }m\ge k;
\end{eqnarray}
and ({\em $h$-adic weak commutativity})
for any $a,b,c\in A[[h]]$ and for any nonnegative integer $n$,
there exists a nonnegative integer $k$ (depending on $a,b$ and $n$) such that
\begin{eqnarray}
(x_{1}-x_{2})^{k}[Y_{h}(a,x_{1}),Y_{h}(b,x_{2})]c
\in h^{n}A[[h]][[x_{1},x_{1}^{-1},x_{2},x_{2}^{-1}]]
\end{eqnarray}
and ({\em $h$-adic weak associativity})
for any $a,b,c\in A[[h]]$ and for any nonnegative integer $n$,
there exists a nonnegative integer $l$ (depending on $a,c$ and $n$) such that
\begin{eqnarray}
(x_{0}+x_{2})^{l}\left( Y_{h}(a,x_{0}+x_{2})Y_{h}(b,x_{2})c
-Y_{h}(Y_{h}(a,x_{0})b,x_{2})c\right)
\in h^{n}A[[h]][[x_{0},x_{0}^{-1},x_{2},x_{2}^{-1}]].
\end{eqnarray}
\ep

{\bf Proof.} Notice that we have
\begin{eqnarray}
\cap_{n\ge 0} h^{n}A[[h]]=0.
\end{eqnarray}
In view of this, for $a\in A[[h]]$, $Y_{h}({\bf 1},x)a=a$ 
if and only if for every positive integer $n$, $Y_{h}(a,x)a-a\in
h^{n}A[[h]][[x,x^{-1}]]$, which is equivalent to
$\bar{Y}_{n}({\bf 1}+h^{n}A[[h]],x)(a+h^{n}A[[h]])=a+h^{n}A[[h]]$.
Similarly, the equivalence on the creation property is clear.
As the Jacobi identity for a vertex algebra (over $\C[[h]]$)
amounts to weak commutativity and weak associativity
(see [DL], [Li2], [LL]), 
the Jacobi identity for vertex algebras $(A[[h]]/h^{n}A[[h]],
\bar{Y}_{n},{\bf 1}+h^{n}A[[h]])$ for $n\ge 1$ amounts to
the $h$-adic weak commutativity and associativity.
$\;\;\;\;\Box$

\br{r-topological-voa}
{\em Notice that in view of (\ref{e-hadic-truncation})
the following Jacobi identity for $Y_{h}$
\begin{eqnarray}
& &x_{0}^{-1}\delta\left(\frac{x_{1}-x_{2}}{x_{0}}\right)
Y_{h}(a,x_{1})Y_{h}(b,x_{2})c
-x_{0}^{-1}\delta\left(\frac{x_{2}-x_{1}}{-x_{0}}\right)
Y_{h}(b,x_{2})Y_{h}(a,x_{1})c\nonumber\\
&=&x_{2}^{-1}\delta\left(\frac{x_{1}-x_{0}}{x_{2}}\right)
Y_{h}(Y_{h}(a,x_{0})b,x_{2})c
\end{eqnarray}
holds in $A[[h]] [[x_{0}^{\pm 1},x_{1}^{\pm 1},x_{2}^{\pm 1}]]$,
where for $m,n,k\in \Z$, the coefficients of 
$x_{0}^{m}x_{1}^{n}x_{2}^{k}$
in the three main terms are infinite convergent sums in $A[[h]]$,
unlike the coefficients of $x_{0}^{m}x_{1}^{n}x_{2}^{k}$ in
the usual (algebraic) Jacobi identity, which are finite sums.}
\er

\br{rbraided-voas}
{\em In [EK], Etingof and Kazhdan introduced a notion of
braided VOA over $\C[[h]]$.
In fact, one can show that $h$-adic vertex algebras are exactly those 
braided VOA's over $\C[[h]]$ with ${\cal{S}}=1$.}
\er

Let $U$ be a vertex algebra over $\C[h]$, or 
an $h$-adic vertex algebra over $\C[[h]]$.
Then $hU$ is an ideal of $U$, so that $U/hU$ is a 
vertex algebra over $\C[h]$ (with $h$ acting as zero). 
We may naturally consider $U/hU$ as a vertex algebra over $\C$.
The following result is due to [FB]:

\bp{pfact}
Let $U$ be either an $h$-adic vertex algebra over $\C[[h]]$ or 
a vertex algebra over $\C[h]$ such that $\ker_{U}h\subset hU$.
Assume that that $U/hU$ is a commutative vertex algebra.
Then $(U/hU, \cdot,Y_{-}, \partial)$ carries the structure of
a vertex Poisson algebra where
\begin{eqnarray}
& &(u+hU)\cdot (v+hU)=u_{-1}v+hU\\
& &Y_{-}(u+hU,x)(v+hU)={1\over h}\Sing Y(u,x)v+ hU[[x,x^{-1}]]\\
& &\partial (u+hU)={\cal{D}}u+hU=u_{-2}{\bf 1}+hU
\end{eqnarray}
for $u,v\in U$. 
\ep

\br{rpfact}
{\em Proposition \ref{pfact} is a variant of Proposition 15.2.4
of [FB] which states that if $V^{\epsilon}$ is a vertex algebra over
$\C[\epsilon]/(\epsilon^{2})$ and a flat module over 
$\C[\epsilon]/(\epsilon^{2})$ such that $V^{0}=V^{\epsilon}/\epsilon
V^{\epsilon}$ is commutative, then $V^{0}$ naturally acquires the
structure of a vertex Poisson algebra. The variant Proposition
\ref{pfact} follows from the same proof of [FB] where the assumption 
$\ker h\subset hU$ guarantees that 
$Y_{-}(u+hU,x)(v+hU)={1\over h}\Sing Y(u,x)v+ hU[[x,x^{-1}]]$ 
is well defined.}
\er

Proposition \ref{pfact} gives another way 
to construct vertex Poisson algebras (over
$\C$) from vertex algebras (over $\C[[h]]$) 
(recall Proposition \ref{pfilteredvoa}).
In the next proposition we show that
the vertex Poisson algebra $\gr_{E}V$ associated 
to a filtered vertex algebra $(V,E)$ can be realized as
a vertex Poisson algebra $U/hU$ for some vertex algebra $U$ 
over $\C[h]$.

The following result is classical in nature:

\bp{pdeformationfilteredva}
Let $(V,E)$ be a filtered vertex algebra over $\C$.
For $n\ge 0$, let $U_{n}$ be a subspace of $E_{n}$ such that
$E_{n}=U_{n}\oplus E_{n-1}$, so that
\begin{eqnarray}
V=U_{0}\oplus U_{1}\oplus \cdots =\coprod_{n\ge 0}U_{n}.
\end{eqnarray}
For $n\ge 0$, denote by $p_{n}$ the projection map of $V$ onto $U_{n}$.
Set
\begin{eqnarray}
V_{E}[h]=V[h]=\C [h]\otimes V,
\end{eqnarray}
as a $\C[h]$-module, which is free. 
For $u\in U_{m},\;v\in U_{n}$, we define
\begin{eqnarray}
Y_{h}(u,x)v=\sum _{j=0}^{m+n}h^{m+n-j}p_{j}(Y(u,x)v)
=\sum_{r\in \Z}\sum _{j=0}^{m+n}h^{m+n-j}p_{j}(u_{r}v)x^{-r-1},
\end{eqnarray}
and then extend the definition $\C[h]$-linearly to $V_{E}[h]$.
Then $(V_{E}[h],Y_{h}, {\bf 1})$ carries the structure of a vertex algebra 
over $\C[h]$ such that $V_{E}[h]/hV_{E}[h]$ is commutative.
Furthermore, the linear map
\begin{eqnarray}
\psi: V&\rightarrow& V_{E}[h]/(h-1)V_{E}[h]\nonumber\\
v&\mapsto& v+(h-1)V_{E}[h]
\end{eqnarray}
is a vertex algebra isomorphism over $\C$ and the linear map
\begin{eqnarray}
\phi: \gr_{E}V=\coprod_{m\ge 0}E_{m}/E_{m-1}
&\rightarrow& V_{E}[h]/hV_{E}[h]\nonumber\\
v+E_{m-1}\in E_{m}/E_{m-1}&\mapsto& p_{m}(v)+hV_{E}[h]
\end{eqnarray}
is a vertex Poisson algebra isomorphism over $\C$.
\ep

{\bf Proof.} First, for every $n\ge 0$, we have 
$$E_{n}=U_{0}\oplus \cdots \oplus U_{n},$$ 
so that 
\begin{eqnarray}
v=\sum_{k\ge 0}p_{k}(v)
=\sum_{k=0}^{n}p_{k}(v)\;\;\;\mbox{ for }v\in E_{n}.
\end{eqnarray}
Second, in view of properties (\ref{enegativemode})
and (\ref{enonnegativemode}) we have 
\begin{eqnarray}\label{eproof-second}
Y(u,x)v\in E_{m+n}((x))\;\;\;\mbox{ for }u\in E_{m},\; v\in E_{n}.
\end{eqnarray}
Third, $\C[h]\otimes V$ is naturally a 
vertex algebra over $\C[h]$. 

Let $\theta$ be the $\C[h]$-linear endomorphism of $V[h]=\C[h]\otimes V$
defined by
\begin{eqnarray}
\theta(f(h)v)=\sum_{n\ge 0}f(h)h^{n}p_{n}(v)\;\;\;\mbox{ for }f(h)\in
\C[h],\; v\in V.
\end{eqnarray}
In particular,
\begin{eqnarray}
\theta(u)=h^{n}u\;\;\;\mbox{ for }u\in U_{n},\;n\ge 0.
\end{eqnarray}
It is easy to see that $\theta$ is injective.
The image of $\theta$, which is $\sum_{n\ge 0}\C[h]h^{n}U_{n}$,
is a vertex subalgebra (over $\C[h]$) because
for $u\in U_{m},\; v\in U_{n}$,
\begin{eqnarray}
& &Y(h^{m}u,x)(h^{n}v)=h^{m+n}Y(u,x)v
=\sum_{k=0}^{m+n}h^{m+n}p_{k}(Y(u,x)v)\nonumber\\
& &\hspace{2cm}=\sum_{k=0}^{m+n}h^{m+n-k}h^{k}p_{k}(Y(u,x)v).
\end{eqnarray}

Define a $\C[h]$-linear map $Y_{h}'$ from $V[h]$ to
$(\End_{\C[h]}V[h])[[x,x^{-1}]]$ by
\begin{eqnarray}
Y_{h}'(u,x)v=\theta^{-1}Y(\theta(u),x)\theta(v)
\;\;\;\mbox{ for }u,v\in V[h].
\end{eqnarray}
Then $(V[h], Y_{h}',{\bf 1})$ carries the structure of a vertex
algebra over $\C[h]$, which is transported from the vertex
algebra $\theta(V[h])$ through the map $\theta$.
Furthermore, for $u\in U_{m},\; v\in U_{n}$, we have
\begin{eqnarray}
Y_{h}'(u,x)v&=&\theta^{-1}Y(\theta(u),x)\theta(v)
=\theta^{-1}\left(Y(h^{m}u,x)(h^{n}v)\right)\nonumber\\
&=&\sum_{k=0}^{m+n}
\theta^{-1}\left(h^{m+n}p_{k}(Y(u,x)v)\right)\nonumber\\
&=&\sum_{k=0}^{m+n}
h^{m+n-k}p_{k}(Y(u,x)v)\nonumber\\
&=&Y_{h}(u,x)v.
\end{eqnarray}
This shows that $Y_{h}=Y_{h}'$. We also have $\theta ({\bf 1})={\bf
1}$ because ${\bf 1}\in E_{0}=U_{0}$. Therefore,
$(V_{E}[h],Y_{h},{\bf 1})$ 
carries the structure of a vertex algebra over $\C[h]$.
With $(h-1)V_{E}[h]$ being an ideal of $V_{E}[h]$, 
$V_{E}[h]/(h-1)V_{E}[h]$ is naturally a vertex algebra over
$\C[h]$. Clearly,
the linear map $\psi$ from $V$ to $V_{E}[h]/(h-1)V_{E}[h]$
is a linear isomorphism sending ${\bf 1}$ to ${\bf 1}+(h-1)V_{E}[h]$.
Furthermore, for $\in U_{m},\; v\in U_{n}$, we have
\begin{eqnarray}
Y(u,x)v=\sum_{k=0}^{m+n}p_{k}(Y(u,x)v)=Y_{h}(u,x)v|_{h=1}.
\end{eqnarray}
Thus $\psi$ is a vertex algebra isomorphism.

For $\in U_{m},\; v\in U_{n}$, since $u_{r}v\in E_{m+n-1}$ 
for $r\ge 0$, we have
$p_{m+n}(u_{r}v)=0$ for $r\ge 0$, so that
\begin{eqnarray}
Y_{h}(u,x)v=\sum_{r<0}\sum_{k=0}^{m+n}
h^{m+n-k}p_{k}(u_{r}v) x^{-r-1}+h\sum_{r\ge 0}\sum_{k=0}^{m+n-1}
h^{m+n-k-1}p_{k}(u_{r}v) x^{-r-1}.
\end{eqnarray}
Thus $Y_{h}(u,x)v\in V_{h}[[x]]+ hV_{h}((x))$. It follows that
the quotient vertex algebra $V_{E}[h]/hV_{E}[h]$ is commutative.
Recall that $E_{m}=U_{m}\oplus E_{m-1}$ for $m\ge 0$.
For $u\in U_{m}=E_{m}/E_{m-1},\;
v\in U_{n}=E_{n}/E_{n-1}$, we have
\begin{eqnarray}
Y_{-}(u+E_{m-1},x)(v+E_{n-1})&=&\sum_{r\ge 0}(u_{r}v+E_{m+n-2})x^{-r-1}
\nonumber\\
&=&\sum_{r\ge 0}(p_{m+n-1}(u_{r}v)+E_{m+n-2}) x^{-r-1}
\end{eqnarray}
and
\begin{eqnarray}
Y_{-}(u+hV_{E}[h],x)(v+hV_{E}[h])
&=&{1\over h}\Sing Y_{h}(u,x)v +hV_{E}[h]\nonumber\\
&=&\sum_{r\ge 0}\sum_{k=0}^{m+n-1}
h^{m+n-k-1}\left(p_{k}(u_{r}v)+hV_{E}[h]\right) x^{-r-1}\nonumber\\
&=&\sum_{r\ge 0}\left(p_{m+n-1}(u_{r}v)+hV_{E}[h]\right) x^{-r-1}.
\end{eqnarray}
We also have
\begin{eqnarray}
(u+E_{m-1})\cdot
(v+E_{n-1})=u_{-1}v+E_{m+n-1}=p_{m+n}(u_{-1}v)+E_{m+n-1}
\end{eqnarray}
and
\begin{eqnarray}
& &(u+hV_{E}[h])\cdot
(v+V_{E}[h])\;\left(=\Res_{x}x^{-1}Y_{h}(u,x)v+hV_{E}[h]\right)\nonumber\\
&=&\sum_{i=0}^{m+n}h^{m+n-k}p_{k}(u_{-1}v)+hV_{E}[h]\nonumber\\
&=&p_{m+n}(u_{-1}v)+hV_{E}[h].
\end{eqnarray}
Now it is clear that $\phi$ is
a vertex Poisson algebra isomorphism (note that ${\bf 1}\in U_{0}=E_{0}$).
$\;\;\;\;\Box$

\br{rdependence}
{\em Technically speaking, $V_{E}[h]$ is not a well defined
function of $(V,E)$ because the vertex algebra $V_{E}[h]$ 
also depends on the decompositions $E_{n}=\coprod_{i=0}^{n}U_{n}$.
However, it is straightforward to check that
different decompositions give rise to isomorphic vertex algebras $V_{E}[h]$.}
\er

\br{rCGconstruction}
{\em Let $(V,E)$ be a filtered vertex algebra. Consider the subspace
\begin{eqnarray}
U=\coprod_{n\in \N}E_{n}h^{n}\subset V[h].
\end{eqnarray}
Since $E_{n}\subset E_{n+k}$ for $n,k\ge 0$, we have
$$E_{n}h^{n+k}\subset E_{n+k}h^{n+k}\subset U,$$
so that $U$ is indeed a $\C[h]$-subspace of $V[h]$.
In view of (\ref{eproof-second}), $U$ is a vertex subalgebra of $V[h]$.
Clearly, $\ker h=0\subset hU$.
For $u\in E_{m},\; v\in E_{n},\; m,n,i\ge 0$, we have
\begin{eqnarray}
h^{m+n}u_{i}v\in h^{m+n}E_{m+n-1}=h (h^{m+n-1}E_{m+n-1}),
\end{eqnarray}
Thus $U/hU$ is commutative. 
We have 
\begin{eqnarray}
U/hU=\coprod_{n\ge 0}E_{n}h^{n}/E_{n-1}h^{n}
= \coprod_{n\ge 0}E_{n}/E_{n-1}=\gr_{E}V,
\end{eqnarray}
as a vector space over $\C$. 
(This is exactly the classical construction (see [CG]).)
Similar to what we did in the proof of Proposition
\ref{pdeformationfilteredva}, one can show 
 that $U/hU=\gr_{E}V$ as a
vertex Poisson algebra and one can also show that 
the evaluation map with $h=1$ from $V[h]$ to $V$
gives rise to a vertex algebra
isomorphism from $U/(h-1)U$ onto $V$.
But, to identify $U$ with $V[h]$ one will have to
use uncanonical decompositions $E_{n}=E_{n-1}\oplus U_{n}$ for $n\ge 0$
just as in Proposition \ref{pdeformationfilteredva}.}
\er

Let $(R,Y_{-},\partial)$ be a vertex Lie algebra 
over $\C$. Extend $\partial$ to a $\C[h]$-linear endomorphism of
$R[h]$ and extend $Y_{-}$ to a $\C[h]$-bilinear map on $R[h]\otimes
R[h]$. Then $(R[h],Y_{-},\partial)$ is a vertex Lie
algebra over $\C[h]$.
It is straightforward to see that
$(R[h], hY_{-},\partial)$ is a vertex Lie algebra 
over $\C[h]$, just as with Lie algebras. 
Let ${\cal{V}}_{h}(R)$ be the corresponding
vertex algebra over $\C[h]$ associated to the vertex Lie algebra 
$(R[h], hY_{-},\partial)$ and let $Y_{h}$ denote 
the vertex operator map:
\begin{eqnarray}
Y_{h}: {\cal{V}}_{h}(R)\rightarrow 
\Hom _{\C[h]}({\cal{V}}_{h}(R), {\cal{V}}_{h}(R)((x))).
\end{eqnarray}
As with ${\cal{V}}(R)$, $R$ is a $\C$-subspace of ${\cal{V}}_{h}(R)$.
So we have natural $\C$-linear maps from $R$ to any quotient spaces 
of ${\cal{V}}_{h}(R)$.
Now we have:

\bt{thvertexliealgebra}
For any vertex Lie algebra $(R,Y_{-},\partial)$ (over $\C$), 
there is a unique vertex algebra isomorphism from ${\cal{V}}(R)$ 
onto the quotient 
vertex algebra ${\cal{V}}_{h}(R)/(h-1){\cal{V}}_{h}(R)$, extending
the natural map from $R$ to ${\cal{V}}_{h}(R)/(h-1){\cal{V}}_{h}(R)$.
The quotient 
vertex algebra ${\cal{V}}_{h}(R)/h{\cal{V}}_{h}(R)$ is commutative
and there is a unique vertex Poisson algebra isomorphism from
$S(R)$ onto the vertex Poisson algebra
${\cal{V}}_{h}(R)/h{\cal{V}}_{h}(R)$ obtained in Proposition
\ref{pfact}, extending the natural map from $R$ to 
${\cal{V}}_{h}(R)/h{\cal{V}}_{h}(R)$.
\et

{\bf Proof.} For $a,b\in R\subset R[h]$, we have
\begin{eqnarray}
& &\Sing Y_{h}(a,x)b=hY_{-}(a,x)b=Y_{-}(a,x)b+(h-1)Y_{-}(a,x)b\\
& & \partial a={\cal{D}}a,
\end{eqnarray}
where ${\cal{D}}$ denotes the ${\cal{D}}$-operator of
the vertex algebra ${\cal{V}}_{h}(R)$.
Thus, the $\C$-linear map $f: a\mapsto a+(h-1){\cal{V}}_{h}(R)$ from $R$ to
${\cal{V}}_{h}(R)/(h-1){\cal{V}}_{h}(R)$ is a vertex Lie algebra
homomorphism. In view of Primc's result (Theorem \ref{tprimc}), 
$f$ uniquely extends to
a vertex algebra homomorphism $\bar{f}$ from ${\cal{V}}(R)$ to 
${\cal{V}}_{h}(R)/(h-1){\cal{V}}_{h}(R)$. 
It follows from the Poincar\'e-Birkhoff-Witt theorem
that $\bar{f}$ is an isomorphism.

For $a,b\in R\subset R[h]$, we have
\begin{eqnarray}\label{eyh-formula}
[ Y_{h}(a,x_{1}),Y_{h}(b,x_{2})]
&=&\Res_{x_{0}}x_{2}^{-1}\delta\left(\frac{x_{1}-x_{0}}{x_{2}}\right)
Y_{h}(Y_{h}(a,x_{0})b,x_{2})\nonumber\\
&=&\Res_{x_{0}}x_{2}^{-1}\delta\left(\frac{x_{1}-x_{0}}{x_{2}}\right)
Y_{h}(\Sing_{x_{0}} Y_{h}(a,x_{0})b,x_{2})\nonumber\\
&=&\Res_{x_{0}}x_{2}^{-1}\delta\left(\frac{x_{1}-x_{0}}{x_{2}}\right)
hY_{h}(Y_{-}(a,x_{0})b,x_{2}).
\end{eqnarray}
Since $R$ is a generating space of ${\cal{V}}_{h}(R)$, 
$\{a+h{\cal{V}}_{h}(R)\;|\;a\in R\}$ is a 
generating subspace of ${\cal{V}}_{h}(R)/h{\cal{V}}_{h}(R)$. 
Notice that
if $S$ is a generating subset of a vertex algebra $V$ such that
$[Y(u,x_{1}),Y(v,x_{2})]=0$ for $u,v\in S$, then $V$ is commutative.
In view of this and (\ref{eyh-formula}),
${\cal{V}}_{h}(R)/h{\cal{V}}_{h}(R)$ is commutative.
In view of Proposition \ref{pfact}, ${\cal{V}}_{h}(R)/h{\cal{V}}_{h}(R)$
is naturally a vertex Poisson algebra. Since
\begin{eqnarray}
{1\over h}\Sing Y_{h}(a,x)b=Y_{-}(a,x)b\;\;\mbox{ for }a,b\in R,
\end{eqnarray}
the map $g: a\mapsto a+h{\cal{V}}_{h}(R)$ from $R$ to 
${\cal{V}}_{h}(R)/h{\cal{V}}_{h}(R)$ is a vertex Lie algebra 
homomorphism. It follows from Proposition \ref{puniversal-pva}
that $g$ uniquely extends to a
vertex Poisson algebra homomorphism $\bar{g}$ from $S(R)$ into
${\cal{V}}_{h}(R)/h{\cal{V}}_{h}(R)$.
Again, it follows from the Poincar\'e-Birkhoff-Witt theorem
that $\bar{g}$ is an isomorphism.
$\;\;\;\;\Box$

Using the same arguments (with Remark \ref{r-on-primc} in place of 
Theorem \ref{tprimc}) we have the following generalization:

\bt{tvpa-form-construction-2}
Let $(R,Y_{-},\partial)$ be a vertex Lie algebra and let
$\lambda$ be a partially defined linear functional on $\ker_{R} \partial$.
Define ${\cal{V}}_{h}^{\lambda}(R)$ to be the quotient vertex 
algebra of ${\cal{V}}_{h}(R)$ modulo the ideal generated by
$a-\lambda(a)$ for $a\in D_{\lambda}\subset\ker_{R}\partial\subset R[h]$. Then
there exists  a (unique) vertex algebra isomorphism map
$f_{\lambda}$ from
${\cal{V}}_{\lambda}(R)$ to
${\cal{V}}_{h}^{\lambda}(R)/(h-1){\cal{V}}_{h}^{\lambda}(R)$ 
such that 
\begin{eqnarray}
f_{\lambda}\psi_{\lambda}(a)=a+(h-1){\cal{V}}_{h}^{\lambda}(R)
\;\;\;\mbox{ for }a\in R,
\end{eqnarray}
where $\psi_{\lambda}$ is the quotient map from ${\cal{V}}(R)$ to
${\cal{V}}_{\lambda}(R)$.
There exists a (unique) vertex Poisson algebra isomorphism 
map $g_{\lambda}$ from $S_{\lambda}(R)$ to 
${\cal{V}}^{\lambda}(R)/h{\cal{V}}_{h}^{\lambda}(R)$ such that
\begin{eqnarray}
g_{\lambda}\phi_{\lambda}(a)=a+h{\cal{V}}_{h}^{\lambda}(R)
\;\;\;\mbox{ for }a\in R,
\end{eqnarray}
where $\phi_{\lambda}$ is the quotient map from $S(R)$ to
$S_{\lambda}(R)$.$\;\;\;\;\Box$
\et

\br{rFB-theorem}
{\em Let $R=\C[\partial]\otimes \g\oplus \C c$ be the
vertex Lie algebra associated with a Lie algebra $\g$ equipped with 
a symmetric invariant bilinear form $\<\cdot,\cdot\>$ 
(see Example \ref{raffinenonzeroc}).
We have $\ker_{R}\partial =\C c$.
Let $\ell$ be any complex number and let $\lambda_{\ell}$ be the linear
functional on $\C c$ such that $\lambda_{\lambda}(c)=\ell$.
Then ${\cal{V}}_{\lambda_{\ell}}(R)$ is isomorphic to
$V_{\hat{\g}}(\ell,0)$ as in [LL] and to $V_{\ell}(\g)$ as in [FB].
In this special case, a version of Theorem \ref{tvpa-form-construction-2}
was obtained in [FB] (Theorem 15.3.2).}
\er

Motivated by Proposition \ref{pdeformationfilteredva}
we introduce the following notion:

\bd{ddeformationquantization}
{\em Let $A$ be a vertex Poisson algebra and let $V$ be a vertex algebra
over $\C$. We say that $V$ is a {\em deformation of $A$}
if there exists a vertex algebra $P[h]$ over $\C[h]$ such that
$P[h]/hP[h]$ is commutative,
$P[h]/hP[h]\simeq A$ and $P[h]/(h-1)P[h]\simeq V$.}
\ed

It follows from Proposition \ref{pdeformationfilteredva} that
for any filtered vertex algebra $(V,E)$, $V$ is 
a deformation of the associated vertex Poisson algebra $\gr_{E}V$.
In view of Theorems \ref{thvertexliealgebra} and 
\ref{tvpa-form-construction-2},
vertex algebras ${\cal{V}}(R)$ and ${\cal{V}}_{\lambda}(R)$ are
deformations of vertex Poisson algebras $S(R)$ and $S_{\lambda}(R)$,
respectively. 

The notion of formal deformation of an algebra (see [G])
naturally suggests the following notion of formal deformation
of a vertex algebra:

\bd{d*-product-on-va}
{\em Let $(V,Y,{\bf 1})$ be a vertex algebra. A {\em formal deformation}
of $V$ is a $\C$-bilinear map $Y_{h}$ from $V\otimes V$
to $V[[h]][[x,x^{-1}]]$ such that
$(V[[h]], Y_{h},{\bf 1})$ carries the structure 
of an $h$-adic vertex algebra over $\C[[h]]$, where $Y_{h}$ extends 
canonically to a continuous $\C[[h]]$-bilinear map 
on $V[[h]]\otimes V[[h]]$, such that for $u,v\in V$,
\begin{eqnarray}
Y_{h}(u,x)v=Y(u,x)v+hY_{1}(u,x)v+h^{2}Y_{2}(u,x)v+\cdots,
\end{eqnarray}
where $Y_{i}$ are linear maps from $V$ to $\Hom (V,V((x)))$.}
\ed

One always has the $h$-adic vertex algebra structure over $V[[h]]$ with
$Y_{h}$ being the natural extension of the vertex operator map $Y$ on $V$.
This deformation is called the {\em trivial} deformation.
Two formal deformations $Y_{h}^{(1)}$ and $Y_{h}^{(2)}$ on $V$ 
are said to be {\em equivalent} if 
there exists a continuous $\C[[h]]$-linear isomorphism $F$ of $V[[h]]$ such that
$F({\bf 1})={\bf 1}$, $F(v)\in v+hV[[h]]$ for $v\in V$ and such that
\begin{eqnarray}
F(Y_{h}^{(1)}(u,x)v)=Y_{h}^{(2)}(F(u),x)F(v)\;\;\;\mbox{ for }u,v\in V.
\end{eqnarray}
In particular, $F$ is an $h$-adic-vertex-algebra isomorphism from
$(V[[h]],Y_{h}^{(1)}, {\bf 1})$ 
to $(V[[h]],Y_{h}^{(2)},{\bf 1})$. 
The  $\C[[h]]$-linear isomorphism $F$
necessarily has the form
\begin{eqnarray}
\Psi= 1+ hf_{1}+h^{2}f_{2}+\cdots
\end{eqnarray}
with $f_{n}\in \End_{\C}V$ for $n\ge 0$.
A vertex algebra $V$ is said to be {\em rigid} (in the category of
vertex algebras) if any formal deformation is equivalent to the
trivial deformation.

\br{rconjecture}
{\em In the classical case, rigidity is closely related to
the complete reducibility of certain modules. In [Z1-2] (cf. [FZ]), 
a notion of rationality of a vertex operator algebra
$V$ was introduced in terms of the complete reducibility of 
certain $V$-modules
(cf. [DLM1]). In principle, one should be able to show that 
rational vertex operator algebras in the sense of [Z1-2] 
are rigid in the category of vertex algebras.}
\er

\br{rcommentscommutativeobjects}
{\em In view of certain commutativity and associativity properties
(see [FLM], [FHL], [DL], [Li2]), vertex algebras are analogous 
to commutative associative algebras (cf. [B2]).
Certain noncommutative analogues of the
notion of usual vertex algebra were studied in [B2], [BK] and [Li4],
where usual vertex algebras are the commutative objects
in a certain sense. Quantum vertex operator 
algebras studied in [EK] can be considered certain formal deformations of
vertex operator algebras in the category of those noncommutative
analogues of vertex algebras.}
\er

Let $Y_{h}$ be a formal deformation on a vertex algebra $V$. 
Since ${\bf 1}$ is also the vacuum
vector of the $h$-adic vertex algebra, from 
Proposition \ref{pexplicit-version} we have 
$Y_{h}({\bf 1},x)=1$, that is, $Y_{i}({\bf 1},x)=0$ for $i\ge 1$. We
also have $Y_{h}(v,x){\bf 1}\in V[[h,x]]$. In general,
$Y_{h}(v,x){\bf 1}$ may depend on $h$.
Now we assume that
\begin{eqnarray}
Y_{h}(v,x){\bf 1}\in V[[x]]\;\;\;\mbox{  for }v\in V,
\end{eqnarray}
i.e.,
\begin{eqnarray}
Y_{i}(v,x){\bf 1}=0\;\;\;\mbox{  for }v\in V,\; i\ge 1.
\end{eqnarray}
Then $Y_{h}(v,x){\bf 1}=Y(v,x){\bf 1}=e^{x{\cal{D}}}v$.
Since $Y_{h}(v,x){\bf 1}=e^{x{\cal{D}}_{h}}v$ for $v\in V[[h]]$,
the ${\cal{D}}$-operator ${\cal{D}}_{h}$ of $V[[h]]$ on $V$ agrees with 
the ${\cal{D}}$-operator ${\cal{D}}$ of $V$. Thus, ${\cal{D}}_{h}$
is the natural $\C[[h]]$-linear extension of ${\cal{D}}$.
{}From the skew symmetry of the vertex algebra $V[[h]]$ 
we obtain
$$Y_{n}(u,x)v=e^{x{\cal{D}}}Y_{n}(v,-x)u\;\;\;\mbox{ for }u,v\in V,\;
n\ge 0.$$
Thus we have:

\bp{pvertexstarproduct1}
Let $Y_{h}$ be a formal deformation of a vertex algebra $V$
with the property
\begin{eqnarray}\label{estrongcreation}
Y_{h}(v,x){\bf 1}\in V[[x]]\;\;\;\mbox{  for }v\in V.
\end{eqnarray}
Then 
\begin{eqnarray}
& &Y_{h}({\bf 1},x)v=v,\\
& &Y_{h}(v,x){\bf 1}=e^{x{\cal{D}}}v,\\
& &Y_{n}(u,x)v=e^{x{\cal{D}}}Y_{n}(v,-x)u\;\;\;\mbox{ for }u,v\in V,\;
n\ge 0,
\end{eqnarray}
where ${\cal{D}}$ is the ${\cal{D}}$-operator of 
vertex algebra $V$.$\;\;\;\;\Box$
\ep

\br{rrecall-star-product}
{\em We here recall the notion of $*$-deformation
 (see [BFFLS], [BW]) to compare with the vertex analogue.  
Let $A$ be a unital commutative associative algebra (over $\C$)
and let $*_{h}$ be a family of associative multiplications on $A$
given by a formal power series
\begin{eqnarray}
a*_{h}b=\sum_{j\ge 0}B_{j}(a,b)h^{j}
\end{eqnarray}
where each $B_{j}: A\times A\rightarrow A$ is a bilinear map. 
Then $*_{h}$ is called a $*$-deformation of $A$
if for $a,b,c\in A$,
 
1. $B_{0}(a,b)=ab$ (the product in $A$).

2. $B_{j}(b,a)=(-1)^{j}B_{j}(a,b)$

3. $B_{j}(1,a)=0$ for $j\ge 1$

4. $B_{i}$ are bidifferential
operators (i.e., bilinear maps $A\times A\rightarrow A$ are
differential operators with respect to each argument of globally
bounded order). 

5. $(a*_{h}b)*_{h}c=a*_{h}(b*_{h}c)$.}
\er

In view of Proposition \ref{pvertexstarproduct1} and
Remark \ref{rrecall-star-product}, we call a formal deformation 
$Y_{h}$ of a vertex algebra $V$ with the property (\ref{estrongcreation})
a {\em $*$-deformation} of $V$.
Furthermore, motivated by
Frenkel and Ben-Zvi's proposition (Proposition \ref{pfact}) we define
the following notion:

\bd{dstarproduct}
{\em Let $(A,Y_{-},\partial)$ be a vertex Poisson algebra (over $\C$). 
A {\em $*$-deformation  of $A$} is an $h$-adic vertex algebra structure 
$(Y_{h}, {\bf 1})$ on $A[[h]]$ with ${\bf 1}=1$ such that
$Y_{h}(a,x)1\in A[[x]]$ for $a\in A$ and such that for $a,b\in A$,
\begin{eqnarray}
Y_{h}(a,x)b=(e^{x\partial}a)b+h Y_{1}(a,x)b+ O(h^{2}),
\end{eqnarray}
where
\begin{eqnarray}
Y_{-}(a,x)b=\Sing_{x} Y_{1}(a,x)b\;\;\;\mbox{ for }a,b\in A.
\end{eqnarray}}
\ed

\br{rconsequenceofdefinition}
{\em Let $Y_{h}$ be a $*$-deformation of a vertex Poisson algebra 
$(A,Y_{-},\partial)$. Then
\begin{eqnarray}
& &(e^{x\partial}a)b=\lim_{h\rightarrow 0}Y_{h}(a,x)b\\
& &Y_{-}(a,x)b=\Sing_{x} Y_{1}(a,x)b
=\lim_{h\rightarrow 0}{1\over h}\Sing Y_{h}(a,x)b
\;\;\;\mbox{ for }a,b\in A.
\end{eqnarray}
Furthermore, we have
\begin{eqnarray}
{1\over h}[Y_{h}(a,x_{1}),Y_{h}(b,x_{2})]
&=&\Res_{x_{2}}x_{0}^{-1}\delta\left(\frac{x_{1}-x_{0}}{x_{2}}\right)
{1\over h}Y_{h}(Y_{h}(a,x_{0})b,x_{2})\nonumber\\
&=& \Res_{x_{2}}x_{0}^{-1}\delta\left(\frac{x_{1}-x_{0}}{x_{2}}\right)
{1\over h}Y_{h}(\Sing Y_{h}(a,x_{0})b,x_{2}).
\end{eqnarray}}
\er

As with Poisson algebras,
the fundamental problem is about
the existence, uniqueness and construction of 
$*$-deformation of each vertex Poisson algebra.

\end{document}